\documentclass[a4paper,prepint]{elsarticle}

\makeatletter
\def\ps@pprintTitle{%
	\let\@oddhead\@empty
	\let\@evenhead\@empty
	\let\@oddfoot\@empty
	\let\@evenfoot\@oddfoot
}
\makeatother

\usepackage{amsfonts,amsmath,amsthm,amstext,amssymb}
\usepackage{graphicx,color}
\usepackage{enumerate}
\usepackage{eurosym} 
\usepackage{lscape}

\usepackage[usestackEOL]{stackengine}

\usepackage{subcaption}
\usepackage{natbib}

\usepackage{multirow} 
\usepackage{rotating}

\usepackage{marginnote}
\allowdisplaybreaks

\usepackage{ulem}

\newcommand{\bra}[1]{\left[#1\right]} 
\newcommand{\prn}[1]{\left(#1\right)} 
\newcommand{\cur}[1]{\left\{#1\right\}} 
\newcommand{\abs}[1]{\left|#1\right|} 
\newcommand{\ov}[1]{\overline{#1}} 
\newcommand{\mb}[1]{\mathbf{#1}} 

\newcommand{\field}[1]{\mathbb{#1}} 
 
\newcommand{\mc}[1]{\mathcal{#1}}

\newcommand{\st}{\rule{0pt}{2.5ex}}
\newcommand{\str}{\rule{0pt}{3ex}}

\newcommand{\cm}[1]{\textcolor{black} { #1}}

\newcommand{\jp}[1]{\textcolor{blue}{#1}}

\newtheorem{proposition}{Proposition}
\renewenvironment{proof}{\textbf{Proof} \nopagebreak }{\hspace{\stretch{1}}\textbf{q.e.d.} \nolinebreak \rule{0.5em}{0.5em}}

\begin{document}

\begin{frontmatter}
\underline{\textbf{Preprint version submitted to Elsevier \hfill November 2, 2021}}\\
\begin{center}
	{\textbf{\Large Market proliferation and the impact of locational complexity on network restructuring}}
\end{center}	
\begin{center}
	Jes\'us Mar\'ia Pinar-P\'erez, Diego Ruiz-Hern\'andez,  and Mozart B.C. Menezes
\end{center}
\vspace{5mm}
\begin{center}
	{\large Published in \textbf{Applied Mathematical Modelling} (ELSEVIER), December 2021.\\}
\end{center}
\vspace{5mm}
\textbf{Cite as:} Pinar-P\'erez, J. M., Ruiz-Hern\'andez, D., and Menezes, M. B. (2022). Market proliferation and the impact of locational complexity on network restructuring. Applied Mathematical Modelling, 104, 315-338.\\
\begin{center}
	\textbf{DOI:} https://doi.org/10.1016/j.apm.2021.11.031
\end{center}
\vspace{5mm}
\begin{flushleft}
	Article available under the terms of the \textbf{CC-BY-NC-ND} licence
\end{flushleft}

\vfill
\textbf{Preprint version submitted to Elsevier \hfill November 2, 2021}

\title{Market proliferation and the impact of locational complexity on network restructuring}

\author[JePi]{Jes\'us M. Pinar-P\'erez}
\author[DiRu]{Diego Ruiz-Hern\'andez\corref{Diego}}
\author[Menezes]{Mozart B.C. Menezes}

\address[JePi]{CUNEF Universidad, Leonardo Prieto Castro 2, 28040, Madrid, Spain}
\cortext[Diego]{Corresponding author.  \textit{Email address}: d.ruiz-hernandez@sheffield.ac.uk}
\address[DiRu]{Sheffield University Management School, Conduit Road, S10 1FL, Sheffield, UK}
\address[Menezes]{NEOMA Business School, 1 Rue du Mar\'echal Juin, 76130 Mont-Saint-Aignan, France}


\begin{abstract}

 This manuscript investigates the problem of locational complexity, a type of complexity that emanates from a company's territorial strategy.   Using an entropy-based measure for supply chain structural complexity (\textit{pars}-complexity), we develop a theoretical framework for analysing the effects of locational complexity  on the profitability of service/manufacturing networks. 
 The proposed model is used to shed light on the reasons why network restructuring strategies may result ineffective at reducing complexity-related costs.  Our contribution is three-fold. 
 First, we develop a novel mathematical formulation of a facility location problem that integrates the \textit{pars}-complexity measure in the decision process.  
 Second, using this model, we propose a decomposition of the penalties imposed by locational complexity into (a) an intrinsic cost of structural complexity; and (b) an avoidable cost of ignoring such complexity in the decision process.   Such a decomposition is a valuable tool for identifying more effective measures for tackling locational complexity, moreover, it has allowed us to provide an explanation to the so-called \textit{addiction to growth} within the locational context.
 Finally, we propose three alternative strategies that attempt to mimic different approaches used in practice by companies that have engaged in network restructuring processes. The impact of those approaches is evaluated through extensive numerical experiments.
 Our experimental results suggest that network restructuring efforts that are not accompanied by a substantial reduction on the target market of the company, fail at reducing complexity-related costs and, therefore, have a limited impact on the company's profitability.
 
\end{abstract}


\end{frontmatter}

\section{Introduction}

Supply chain complexity is considered a consequence of the proliferation of products, markets, and channels \citep{Mariotti08}, and has been blamed for  eroding firms' capacity for generating profit. \citet{Mariotti08}, observed that while proliferation may lead to an increase in sales, it also brings hidden costs liable to grow faster than revenue. \citet{GeoWil2004} highlighted that the proliferation of products, services and processes in a company, induces large non-value-adding costs --often hidden in overhead structures- that translate in significant profit losses and limit growth. Since then, the problem of supply chain complexity has been addressed by several authors, among them, \citet{Ulr2013},  \citet{ChoSod2014},  \citet{BanBou2016} and \citet{ShaMouAga2016}. More recently \citet{AdaBorBur2020} insisted on the importance of fighting variety-driven complexity as a mechanism for achieving cost-efficiency.

According to \citet{SaeYou1998} complexity can be understood as ``the systemic effect that numerous products, customers, markets... have on activities, overhead structures, and information''. This notion has been narrowed down by \citet{RuiMenAmr2019} for defining supply chain \textit{structural complexity} as the negative effect of the proliferation of products, distribution channels and markets. This complexity stems from strategic choices, and grows as an organisation adds products and/or increases the number of interactions in the supply chain \citep{Hey2007}.

In this manuscript, we address the problem of \textit{locational complexity}, understood as the effect of an increasing number of facilities and their catchment areas (markets) on the operational and financial performance of the company.  \citet{FisGauKle2017} have argued that competitive pressures drive firms to proliferate production, distribution and retail facilities, expanding their network beyond the limits of profitability ``until their chains begin to collapse under their own weight''. The authors refer to this phenomenon as the ``addiction to growth''.

\begin{figure}[!ht]
 \begin{subfigure}{0.5\textwidth}
  \includegraphics[width=6cm]{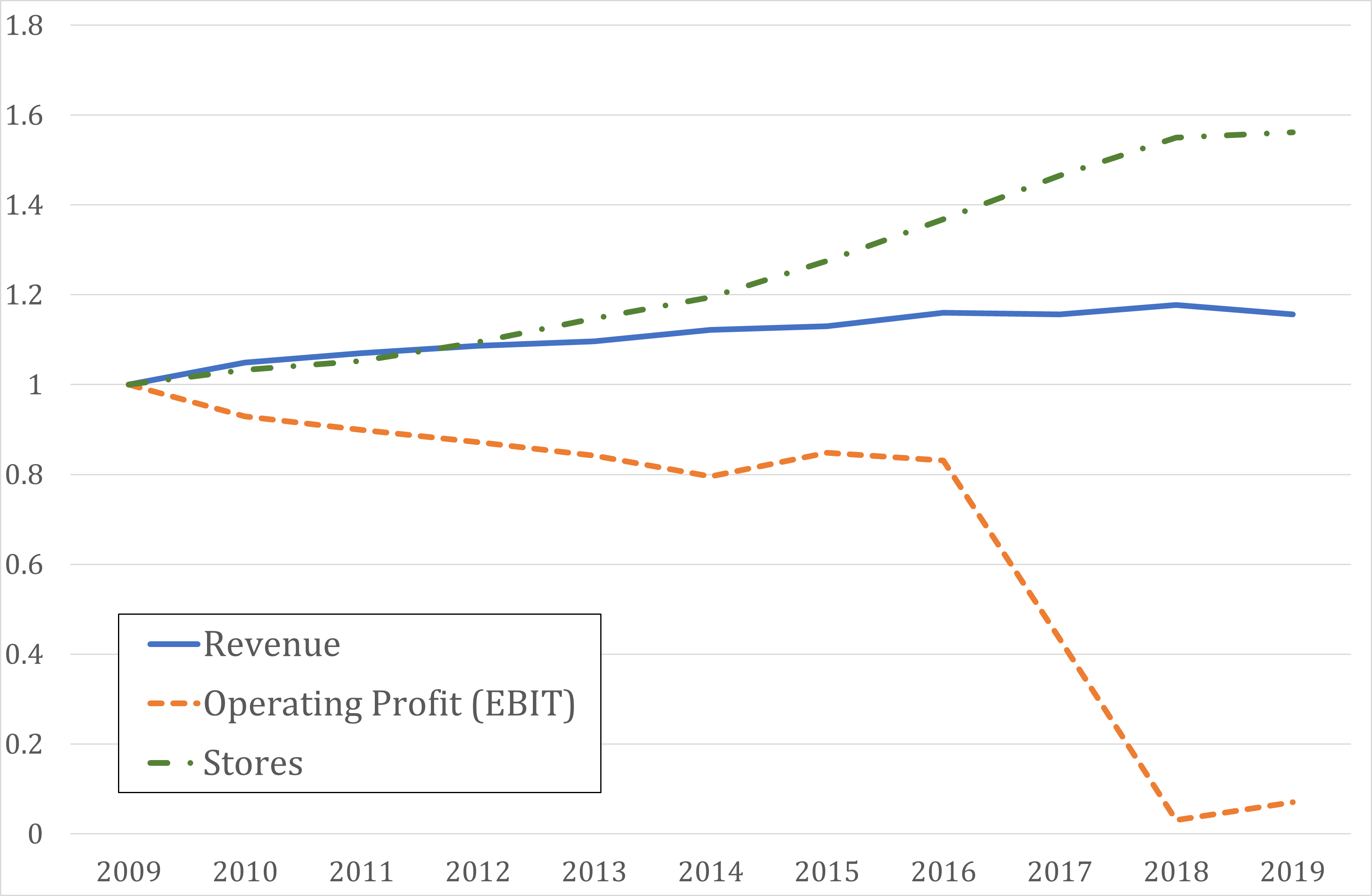}
  \caption{Marks and Spencer (Base 2009)}
 \end{subfigure}
 \begin{subfigure}{0.5\textwidth}
  \includegraphics[width=6cm]{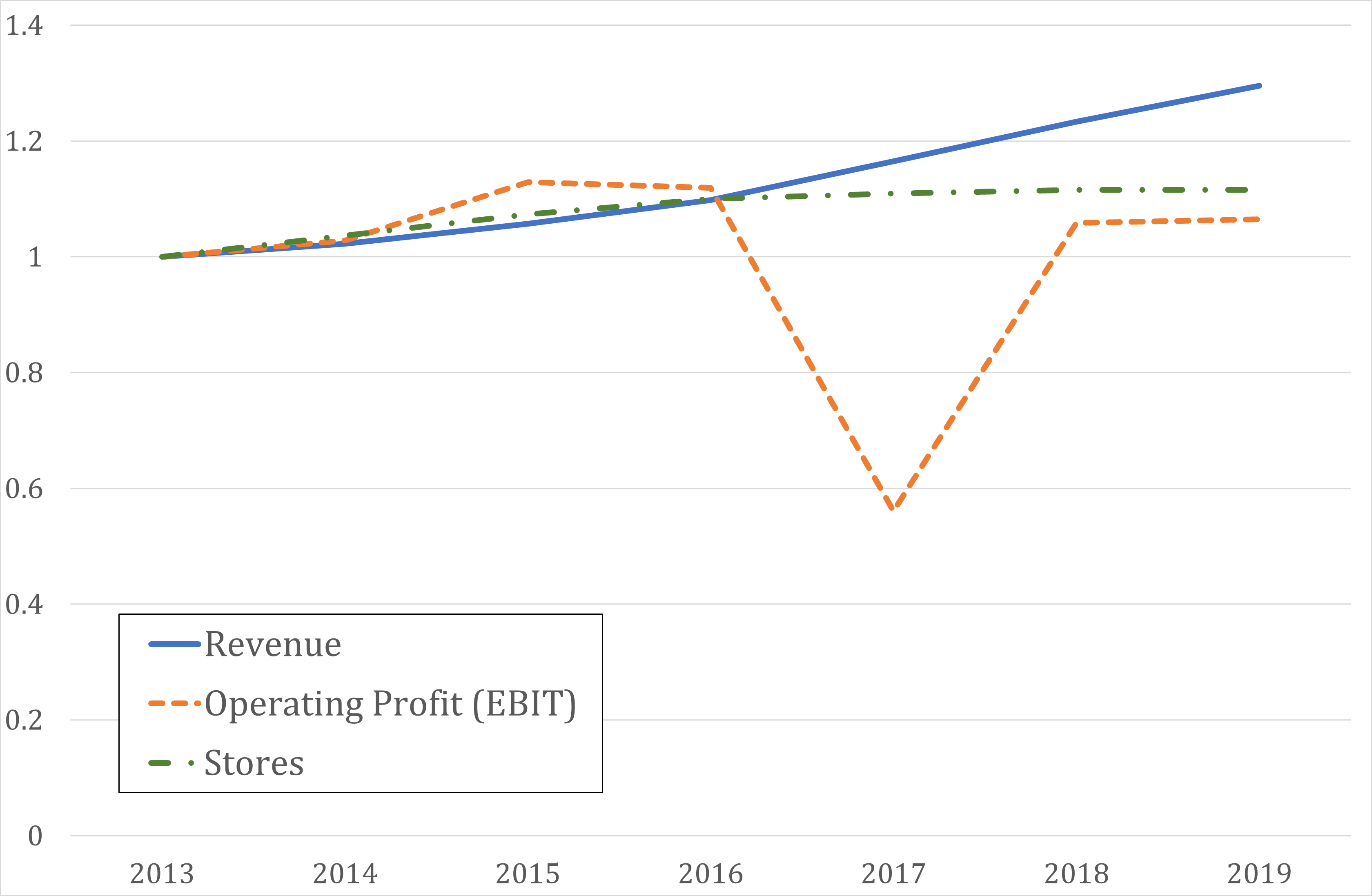}
  \caption{Mercadona (Base 2013)}
 \end{subfigure}
\caption{Operating profit, revenue and stores of two retailers.}\label{fig:retail}
\end{figure}

The negative effect of locational complexity is illustrated in Figure \ref{fig:retail}. The graphs depict the operating profit, revenue and number of stores for two different retailers: Marks \& Spencer (2009-2019) in the UK, and  Mercadona (2013-2019) in Spain. In the first case, despite a sustained growth in the number of stores and in revenue, the operating profit has decreased steadily. In the second case, while the number of branches has increased, and the revenue also shows an increasing pattern, the company's operating profit has remained somehow stagnant over the illustrated period.

It has been argued that the most efficient strategy to mitigate the negative effects of large networks is to eliminate sources of operational and structural complexity that generate non-value-adding costs \citep{ATK2004}. Indeed, \citet{SaeYou1998} found that companies that have managed to reduce variety-driven complexity have used SKU (stock keeping unit) elimination to a degree. Regarding locational complexity, the advice is the same: ``stop opening new stores'' \citep{FisGauKle2017}; and ``stay close to your customer but not too close'' \citep{GeorgeGr2006}. More recenlty, \citet{AitBozGar2016} and \citet{TurAitBox2018}, have highlighted the importance of using appropriate  actions aimed at reducing what they define as \textit{deleterious complexity}, i.e. the elimination of sources of dysfunctional complexity. It is not hard to make a list of companies that have targetted network reduction as a mechanism for surviving in the marketplace. For example, Marks \& Spencer announced the closure of 100 stores by 2022 (BBC News, 03/05/2019). In the UK banking sector, it is expected that 50\% of the total branches will close within the next 10 years (The Guardian, 30/11/2019). In Spain, Santander Bank has revealed plans to close nearly one thousand branches in the short term (Reuters, 13/11/2020).

\begin{figure}[!ht]
 \begin{subfigure}{0.5\textwidth}
  \includegraphics[width=6cm]{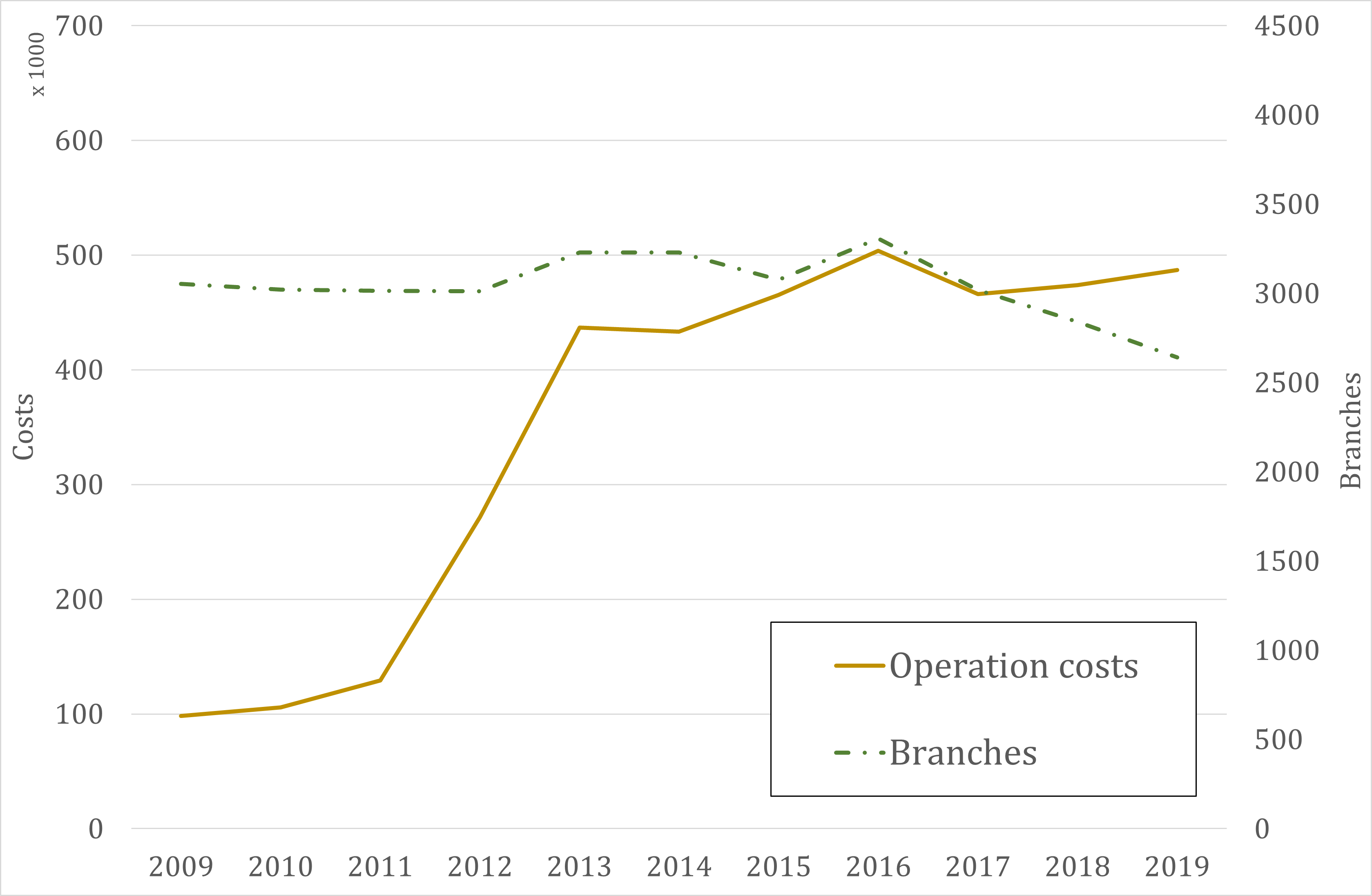}
  \caption{BBVA S.A.}
 \end{subfigure}
 \begin{subfigure}{0.5\textwidth}
  \includegraphics[width=6cm]{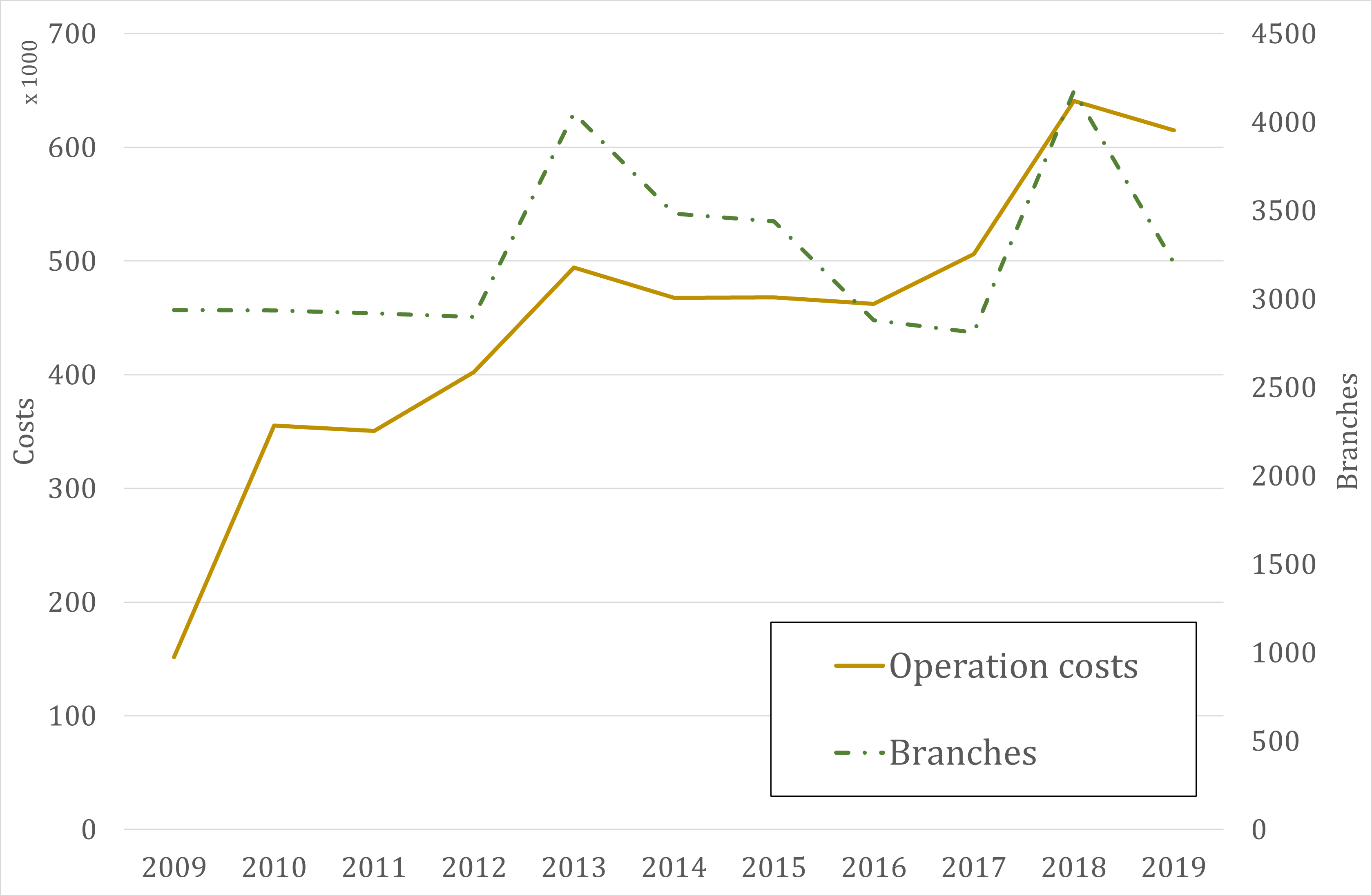}
  \caption{Banco Santander S.A.}
 \end{subfigure}
\caption{Operating costs (left-hand-side axis) and number of branches (right-hand-side axis) of two banking institutions.}
\label{fig:banks}
\end{figure}

However, empirical evidence reveals that network reduction attempts haven't always resulted in considerable cost savings. Figure \ref{fig:banks} shows the operating costs and number of facilities registered over the last 10 years for two Spanish banking groups. Both graphs suggest that the impact of closing branches on the operating expenses has been negligible. Indeed, ignoring the years before 2013, while the number of BBVA branches has been reduced by almost 20\%, the operating costs have continued increasing. With respect to Santander, the sharp reduction in the number of branches observed between 2013 and 2017 (30\% of branches were closed), did not have the expected impact on operating costs\footnote{Please notice that the increase in branches that can be observed for Santander in 2018 is a result of the acquisition of Banco Popular Espa\~nol in 2017. Most of the new acquired branches were closed by 2019 (Source: \texttt{https://www.santander.com/en/press-room/adquisicion-banco-popular}).}.

The aim of this work is to shed light on the forces that hinder the efforts for reducing locational-complexity related cost, and to explore alternative approaches that may help to address this problem. 
The main objective is providing a framework that takes into account complexity-induced costs when assessing the performance of a company's distribution network. 
Based on an entropy-based measure for supply chain structural complexity introduced by \citet{RuiMenAmr2019}, referred to as the \textit{pars}-complexity measure, we propose a novel mechanism for integrating the effect of locational complexity in the decision process. 
The use of this measure has enabled us to propose an explanation for the apparent ineffectiveness of network-restructuring strategies as mechanisms for reducing complexity-related costs.

The rest of the paper is organised as follows: Section \ref{sec:lit} explores the related literature. In Section \ref{sec:form}, we provide the mathematical formulation of the $K$-MedianPlex problem, a $K$-Median-based formulation that introduces the impact of structural complexity in the location problem's objective function. We accompany our formulation with a discussion of the impact of structural complexity on profits. This discussion is further used to suggest an explanation to the problem of addiction to growth as described by \citet{FisGauKle2017}. In Section \ref{sec:ineffective} we present our attempt for explaining the underlying reasons for the ineffectiveness of restructuring processes for reducing costs in closed networks. The methodological approaches aimed at reducing complexity by means of network restructuring are presented in Section \ref{sec:algor}. Section \ref{sec:NumEx} presents the results of our numerical experiments. Finally, Section \ref{sec:conc} concludes the paper.

\section{Related Literature}\label{sec:lit}

Despite the exitence of a vast amount of literature in network restructuring, the locational complexity problem has beem scarcely addressed in the facility location literature. Indeed, most of the available work in supply chain complexity focuses either on the design of alternative measures for complexity (please see Table \ref{tab:comp} for an extensive list of different metrics proposed in literature);  on the quantification of the cost of complexity and its impact on performance \citep{ATK2007a}\nocite{WuFriEfs2007,Sch2008,BozWarFly2009,AdaAlldMue2016,KriLev2018,EkiBay2019}
-\citep{MenRui2020}; or on the analysis of product-variety driven complexity \citep{BleAdb2006a}\nocite{CloJacSwi2008,Chengetal2014,HudMah2019,TraHvaFor2019,ChhHasLau2021}-\citep{MenPin2021}.

Supply chain network complexity has been addressed by \citet{BatPerAll2007}, who use network analysis techniques and an entropy-based measure to quantify the structural complexity of an industrial network. 
\citet{OliElM2018} propose a framework for analysing supply chain robustness and complexity, taking into account the product portfolio structure. \citet{RuiMenAmr2019}, propose an entropy-based mesure of complexity that aims at quantifying structural complexity in terms of the information flows generated by the different echelons in the supply chain.

Regarding the measure of complexity, several authors have proposed alternative metrics. These metrics rely on different types of information, making their implementation in situations or frameworks other than those for which they were built difficult. The most relevant metrics are summarised in Table \ref{tab:comp}, which is a revised version of Table 4 in \cite{RuiMenAmr2019}. Entropy-based measures have also been proposd by \citet{LevPtu2015} and \citet{LevPtu2018}, but with a focus on supply chain risk and have therefore not been included in the table. The more flexible measures seem to be the entropy based ones; among them, we have chosen the \textit{pars}-complexity measure because of its decomposability, ease of calculation, austerity in the use of information, and beacuse it is the only one that seems suitable for measuring locational complexity.

\begin{table}[!ht]
\scriptsize

\begin{tabular}{llp{5.5cm}c}\hline
  \multicolumn{1}{c}{Type of}     &  &                      & \\
 \multicolumn{1}{c}{Complexity} & \multicolumn{1}{c}{Measure} & \multicolumn{1}{c}{Variables} &   Reference \\ \hline
 \str Product & Entropy Based & Segment's share over total sales of the firm &   \cite{Palepu1985} \\
 Manufacturing & Entropy based& Bills of materials, routings, work centres, demand pattern. Periodical data. &  \cite{FriWoo1995} \\ 
 Supply chain & Entropy based &  Imperfectly specified. Forecasted, requested, scheduled and confirmed deliveries. Target and actual production. Materials and information flow. Periodical data.   & \cite{SivEfsShi1999} \\
 Supplier/customer & Entropy based & Imperfectly specified. Flow focused. Demand, production, deliveries. Variations in time and quantities. Periodical data.  &  \cite{Siv2002}\\
 Supply chain & Regression  & Lead time, throughput time, late delivery, tardiness, AMT investment, vertical integration, quality failures, firm size, etc.  &  \cite{VacKlas2002} \\
 Processes & Algebraic & Total value-add time in the process, percentage of defective products, unit processing time, total demand. &  \cite{GeoWil2004} \\
 Supplier/customer & Entropy based & Flow variations (order-forecast, delivery-order, actual-scheduled production), time and/or quantity variations. &   \cite{SivEfsCal2006} \\
 
 Supply chain & Enbtroy-based & Eleven different types of nodes,their interactions, and the monetary flows (costs) among them & \cite{BatPerAll2007}\\
  
 Supply chain & Algebraic & Number of SKUs, markets served, company legal entities, facilities, employees, suppliers, customers. Sales revenues. & \cite{Mariotti08}\\
 
 Supply chain & Regression  & Number of customers, life cycle, number of active material parts, number of products, number of suppliers, percentage of purchases imported, etc. &   \cite{BozWarFly2009} \\
 Supply chain & Entropy based & Expected and actual orders per month (flow based). &   \cite{Isik2010,Isik2011} \\
 
  Supply chain & Statistical  & Expert survey. Number of SKUs, stock locations, employees and years active in market. Several questionnaire topics. &   \cite {LeeGroGoo2013} \\ 
  
  Product   & Algebraic & Number of variants, common elements of components, number of connections. &   \cite{Jac2013} \\
  
  Supply chain   & Entropy based & Percentage contribution to total sales of each combination of SKUs, market, and channel.  &  \cite{RuiMenAmr2019,MenRui2020} \\
  
  Supply chain   & Entropy based & Number of echelons and nodes in the network  &   \cite{LinWanLee2021} \\ \hline
\end{tabular}
\caption{Different complexity measures and their information requirements.}\label{tab:comp}
\end{table}

Given the connection between our study and the literature on facility location and network restructuring, it is worth mentioning certain contributions in these areas that are linked to our research. \citet{AmiBak2017} present a facility location model for closed-loop supply chain network design. 
\citet{CorLopMel2019} propose a multi-stage network design problem for a supply chain with in-house production and outsourcing. \citet{AllYagFat2021} introduce a spatial decision-support methodology to restructure the network of a financial institution. 
 \citet{PouBozYou2021} propose a facility location and network design problem for a healthcare system. Finally, \citet{YavMou2021} expand the available literature on network restructuring by introducing hierarchical facilities in their model.
  
 Unfortunately, the lack of consensus on an operational measure of complexity, has hampered the development of network restructuring models that include the network's locational complexity as a decision variable. The model proposed below constitutes an attempt for filling this gap in literature.

\section{Problem formulation} \label{sec:form}

Discrete and network facility location problems are well known combinatorial problems aimed at finding the best location for a collection of production or distribution facilities with the aim of optimising a given performance measure. One of the most widely used formulations is the $K$-median problem, whose objective is to minimize the cost incurred when serving customers from a set of facilities (see, for example,\citet{Daskin2013}). In this work, we propose a variant of the traditional $K$-median problem that includes a measure of complexity and its related costs, the $K$-MedianPlex problem. The rest of this section is devoted to developing the mathematical formulation of the problem, starting with the presentation of the so called \textit{pars}-complexity measure.

\subsection{\textit{Pars}-Complexity}

\cm{The measure for supply chain structural complexity used in this manuscript is built around a collection $\wp$ of triplets $\cur{SKU,Market,Channel}$ that characterise the flow of products in a firm’s supply chain. Each of these triplets, referred to as \textit{pars} (plural \textit{partes}) has associated a probability value $p_{i} \in \prn{0,1}$ that represents the likelihood that a money unit of revenue comes from the sell by triplet $i\in\wp$. With these elements, we define system's $\wp$ \textit{pars}-complexity as}

\begin{align}\label{eq:Cp}
C_p\prn{\wp}=C_p \left(\mb p\right)  = \sum_{i\in\wp} p_i \log_2\left(\frac{1}{p_i}\right),      
\end{align}
where $\mb p$ is the vector of weights, $p_i$, for $i=1, \ldots, \abs{\wp}$, satisfying $\sum_{i \in \wp} p_i=1$.

\vspace{6pt}

Function $C_p\left(\mb p \right)$ satisfies the following properties:
\begin{enumerate}
	\item $C_p\left(\mb p\right)$ is continuous and concave in $p$.
	\item If $\left|\wp\right|=1$, then the system's complexity is zero, i.e. $C_p\left(\mb p\right)=0$.
	\item $C_p$ attains a maximum when $p_i=\frac{1}{\ \left|\wp\right|}$ for all $i=1,\ ...,\left|\wp\right|$. Such maximum is equal to $\log_2{\left(\left|\wp\right|\right)}$.
\end{enumerate}	
Proof of this properties, together with a discussion around the attributes of the \textit{pars}-complexity measure, can be found in \citet{RuiMenAmr2019}.

The \textit{pars}-complexity metric is simply Shannon's measure of information \citep{Shannon1948}: the expected quantity of information necessary to describe a system's state. As described by \citet{RuiMenAmr2019}, the metric quantifies the information generated by the flow of goods and services in the supply chain, and can thus be used as a proxi for the structural complexity of a supply chain.

\vspace{12pt}

To put the structural complexity problem within a locational complexity framework, we consider a distribution network $\mc S$ where a central manager coordinates a number of regional or local facilities (retail or local distribution centres) serving a set $\mc N$ of demand nodes. For the sake of simplicity and without any loss of generality, we assume $\mc S \subset \mc N$. Each regional facility $\ell$ serves a subset $\mc{N}_\ell$ of demand nodes, where $\mc{N}_\ell= \cur{i \in \mc N: d\prn{i,\ell}\leq \inf_{\ell' \in \mathcal{S}}\cur{d\prn{i,\ell'}}}$ for  $\ell \in \mc S$. Let $\mc{N}^{\mc S}$ represent the collection of allocation sets, $\mc N_\ell$, associated to all facilities $\ell$ in distribution network $\mc S$.
 
Assuming that each demand node $i \in  \mc N$ contributes a proportion $p_i$ of the total demand.  The total complexity of system $\cur{\mc N^{\mc S},\mc S}$ can be written as
\begin{equation}\label{eqn:CompTot}
        C_p\prn{\mc N^{\mc S},\mc S}=\sum_{i \in \mc N} p_{i} \log_2\frac{1}{p_i}=\sum_{\ell \in \mc S} q_\ell \log_2 \frac{1}{q_\ell}+\sum_{\ell \in \mc S} q_{\ell}C_p\prn{\mc N_\ell,\ell}\\
\end{equation}
where $q_\ell=\sum_{i \in \mc N_\ell } p_i$ represents the fraction of the total demand served by facility  $\ell$; and  $C_p\prn{\mc N_\ell,\ell}$ represents the complexity of the single-facility system $\cur{\mc N_\ell,\ell}$. 

The first equality in Equation \eqref{eqn:CompTot} expresses system $\mc{S}$'s \textit{pars}-complexity in terms of the contributions of each demand node. The term after the second equal sign decomposes the total \textit{pars}-complexity in two parts: the first element represents the \textit{pars}-complexity faced by a (hypothetical) central manager, and the second one the weighted average of the \textit{pars}-complexities of the local units.

A consequence of Properties 1 to 3 in a locational context, is that networks with the same number of nodes may have different $C_p$ values depending on the distribution of the demand accross the different nodes. Moreover, smaller networks may show higher complexity values than larger ones.

\vspace{12pt}

\begin{figure}[h!]
\centering
\begin{tabular}{cc}
 \fbox{\includegraphics[width=5.5cm]{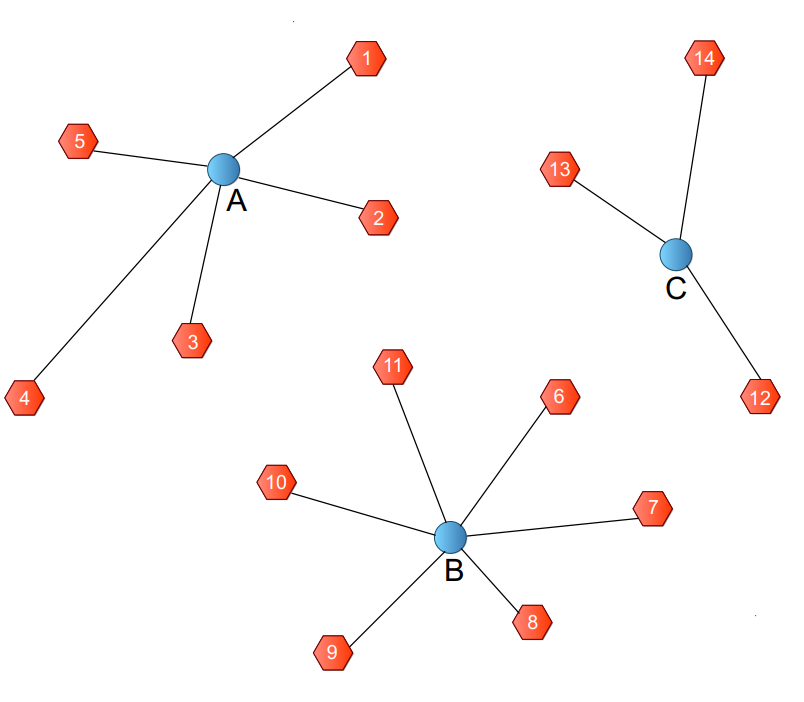}} &
 \fbox{\includegraphics[width=5.5cm]{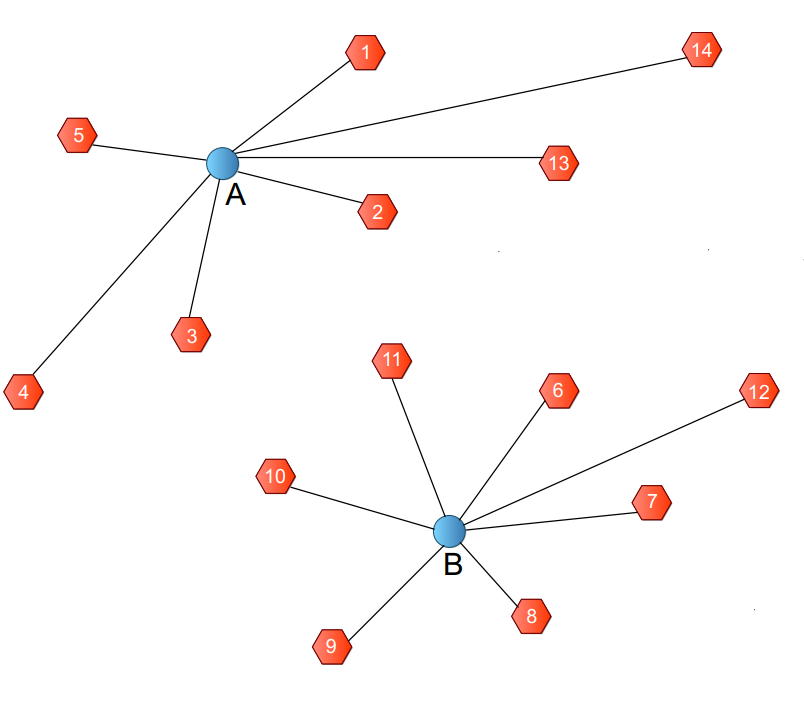}} \\
 \str (a) & (b)
\end{tabular}
\caption{Two alternative configurations of a distribution network with 14 demand nodes and (a) three facilities; or (b) two facilities.} \label{fig:toyex}
\end{figure}

An illustrative example of the use of metric within the locational context is provided in Figure \ref{fig:toyex}. Panel (a) in the figure shows a distribution network consisting of three facilities and fourteen demand nodes. The demand nodes are assigned to their closest server. In this network, $\mc S=\cur{A,B,C}$; $\mc N_A=\cur{1,2,3,4,5}$; $\mc N_B=\cur{6,7,8,9,10,11}$; $\mc N_C=\cur{12,13,14}$; and $\mc N^{\mc S}=\cur{\mc N_A,\mc N_B,\mc N_C}$. The details of the calculation of $C_p\prn{\mc N^{\mc S},\mc S}$ are presented in Table \ref{tab:toyex}. The third column is used to compute the \textit{pars}-complexity of the system according to the left-hand-side term in \eqref{eqn:CompTot}. The fourth column provides the contribution of each demand node to the associated facility, $q^\ell_i=\tfrac{p_i}{q_\ell}, \,\forall  i\in \mc N_\ell, \ell \in \mc S$. Finally, the last column provides the elements for computing $C_p$ as indicated by the right-hand-side term of equation \eqref{eqn:CompTot}.

Panel (b) in Figure \ref{fig:toyex} shows a restructuring of the network presented in panel (a), where facility $C$ is removed and its associated nodes allocated to their second closest facility: $A$ for nodes $13$ and $14$; and $B$ for node $12$. It can be verified that while the total system's complexity remains the same ($C_p=3.624$), the complexities of the remaining facilities  change from $C_p\prn{\mb q^A}=2.202$ and $C_p\prn{\mb q^B}=2.460$ to $2.563$ and $2.697$, respectively. Likewise, the central manager's $C_p$ value (quantity $S1$ in Table \ref{tab:toyex}) changes from $1.485$ to $0.984$.

\begin{table}[h!]
 \begin{tabular}{c|ccccc}\hline
 \str Facility	&	$i$	&	$p_i$	&	$p_i \log_2 \frac{1}{p_i}$	&		$q_i^\ell$	&	$q_i^\ell \log_2 \frac{1}{q_i^\ell}$		\\ \hline
\multirow{7}{*}{A}	&	1	&	0.020	&	0.113	&		0.100	&	0.332		\\
	&	2	&	0.030	&	0.152	&		0.150	&	0.411		\\
	&	3	&	0.070	&	0.269	&		0.350	&	0.530		\\
	&	4	&	0.040	&	0.186	&		0.200	&	0.464		\\
	&	5	&	0.040	&	0.186	&		0.200	&	0.464		\\ \cline{2-6}
\str 	&	\multicolumn{1}{r}{$q_A=$}	&	0.200	&		&	\multicolumn{1}{r}{$C_p(\mb q^A)=$}	&	2.202		\\
\str    &   \multicolumn{1}{r}{$q_A \log_2 \frac{1}{q_A}=$}    &  0.464     &           &  \multicolumn{1}{r}{$q_A\cdot C_p(\mb q^A)=$}     &   0.440           \\ \hline
\multirow{8}{*}{B}	&	6	&	0.090	&	0.313	&		0.180	&	0.445		\\
	&	7	&	0.110	&	0.350	&		0.220	&	0.481		\\
	&	8	&	0.055	&	0.230	&		0.110	&	0.350		\\
	&	9	&	0.030	&	0.152	&		0.060	&	0.244		\\
	&	10	&	0.130	&	0.383	&		0.260	&	0.505		\\
	&	11	&	0.085	&	0.302	&		0.170	&	0.435		\\ \cline{2-6}
\str 	&	\multicolumn{1}{r}{$q_B=$}	&	0.500	&			&	\multicolumn{1}{r}{$C_p(\mb q^B)=$}	&	2.460	\\
\str    &   \multicolumn{1}{r}{$q_B \log_2 \frac{1}{q_B}=$}    &  0.500     &            &  \multicolumn{1}{r}{$q_B\cdot C_p(\mb q^B)=$}     &   1.230      \\ \hline
	\multirow{5}{*}{C}	&	12	&	0.075	&	0.280	&		0.250	&	0.500		\\
	&	13	&	0.111	&	0.352	&		0.370	&	0.531		\\
	&	14	&	0.114	&	0.357	&		0.380	&	0.530		\\ \cline{2-6}
	\str 	&	\multicolumn{1}{r}{$q_C=$}	&	0.300	&		&	\multicolumn{1}{r}{$C_p(\mb q^C)=$}	&	1.561		\\
\str    &   \multicolumn{1}{r}{$q_C \log_2 \frac{1}{q_C}=$}    &  0.521     &     &  \multicolumn{1}{r}{$q_C\cdot C_p(\mb q^C)=$}     &   0.468         \\ \hline \hline
\str	&  $S1:\sum_{\ell }q_\ell \log_2 \frac{1}{q_\ell}=$		&	1.485       &        \multicolumn{2}{r}{$S2:\sum_\ell q_\ell C_p(\mb q^\ell)=$}     &  2.138     \\ \cline{2-6}
\str $C_p(\mb p)$	& \multicolumn{2}{r}{$\sum_{i=1}^{12} p_i \log_2 \frac{1}{p_i}=$}	&	3.624	&			\multicolumn{1}{r}{$S1+S2=$}	&	3.624		\\ \hline
 \end{tabular}
 \caption{Two alternative ways for computing the \textit{pars}-Complexity value of a distribution network with 14 demand nodes and three-facilities.}\label{tab:toyex}
\end{table}

\subsection{Gross profit and the effect of complexity}

It has been argued \citep{SaeYou1998,Mariotti08,FisGauKle2017} that supply chain complexity introduces hidden costs that erode the company's capability for generating profits. In an empirical study, \citet{MenRui2020} found that structural complexity can have an impact in operational profit between 2\% and 10\%  per \textit{pars}-complexity point.  Previously, \citet{ATK2007a} estimated that reducing complexity may lead to an increase of 3 to 5\% in EBIT. In this section we present a mathematical expression for the company's profit that takes into account the effect of complexity.  For completeness, in the following lines we present the main notation used throughout this manuscript.

\vspace{12pt}

\noindent \textbf{Notation}
\begin{description}
\item[$\mc N:$] Set of demand nodes
\item[$\mc S:$] Set of facilities in the distribution network, $\mc S \subset \mc N $  
\item[$\mc{N}_\ell:$] Subset  of demand nodes allocated to facility $\ell \in \mc{S}$ 
\item[$\mc{N}^{\mc S}:$] Collection of all allocation sets $\mc N_\ell \ \forall \ell \in \mc S$, $\mc N^{\mc S} \subseteq \mc N$
\item[$W_i:$] Demand generated by node $i \in \mc N^{\mc S}$
\item[$p_{i}:$] Proportion of total demand generated by node $i \in \mc N^{\mc S}$ 
\item[$\mb p:$] Vector of weights $p_{i}, \, i=1,\ldots, \abs{\mc N^{\mc S}}$ 
\item[$C_p(\mc  N^{\mc S}, \mc  S):$] \textit{pars}-Complexity of system $\cur{\mc N^{\mc S},\mc S}$, also represented by $C_p(\mb p) $
\item[$C_p(\mc  N_\ell, \ell):$] \textit{pars}-Complexity value of the single-facility system $\cur{\mc N_\ell, \ell}$
\item[$d\prn{i,\ell}:$] Distance between the node $i$ and facility $\ell \in \mc S$
\item[$q_\ell:$] Fraction of the total demand served by facility  $\ell \in \mc S$ 
\item[$q_i^\ell:$] Contribution of node $i \in \mc N_\ell$ to the total demand served by facility $\ell \in \mc S$
\item[$\alpha:$] Complexity cost factor 
\item[$\Pi^\circ \prn{\mc N^{\mc S},\mc S}:$] Total profit earned by the distribution system $\cur{\mc N^{\mc S},\mc S}$  
\item[$R \prn{\ell}:$] Gross profits earned by facility $\ell \in \mc S$ 
\item[$C_\alpha \prn{\mc N^{\mc S},\mc S}:$] Cost of complexity for distribution system $\cur{\mc N^{\mc S},\mc S}$
\item[$r:$] Revenue per unit 
\item[$\gamma:$] Transport cost per distance and product unit
\item[$\phi_\ell:$] Fix opening/operating cost of facility $\ell \in \mc S$
\item[$K:$] Number of facilities in system $\mc S$, $K=\abs{\mc S}$
\item[$Z^K_{Plex}:$] Objective function of the $K$-MedianPlex problem
\item[$Z^K:$] Objective function of the $K$-Median problem
\item[$Z^\circ_{Plex}\prn{\mc N^{\mc S},\mc S}:$] Profit function for distribution network $\cur{\mc N^{\mc S},\mc S}$
\end{description}
\vspace{12pt}

Assuming that \textit{pars}-complexity imposes a penalty (or hidden cost) on benefits represented by parameter $0 \leq \alpha <1$, the gross profit earned by a network consisting of $\abs{\mc S}$ facilities is given by

  \begin{align}\label{eq:ProfitLoc}
    \Pi^\circ \prn{\mc N^{\mc S},\mc S}&=  \sum_{\ell \in \mc S} \prn{1-\alpha C_p\prn{\mc N_\ell,\ell}\st}  R \prn{\ell}
   \end{align}
where 
   \begin{align}\label{eq:revenue}
R\left(\ell \right)=\sum_{i\in\mathcal{N}_\ell}{\left(r-\gamma d\prn{i,\ell}\right)\, W_i}, \quad & \ell \in \mathcal{S} 
\end{align}
represents the sum of the profits earned by regional facility $\ell$ when serving all the assigned demand nodes, $i \in \mc N_{\ell}$.

Each demand node's profit is defined as the difference between the revenue per unit, $r$, and the unit transport cost, $\gamma d\prn{i,\ell}$, multiplied by the node's demand, $W_i$. Parameter $\gamma$ represents the transport cost per distance and product unit.  

Expressions \eqref{eq:ProfitLoc} and \eqref{eq:revenue} presume a perfectly competitive market where prices and demand are exogenous (and constant); and the revenue per unit, $r$, is net of production costs. Transport costs are also fixed and distances are invariant.  It is also assumed that parameter $\alpha$ is constant and known, satisfying $1-\alpha C_p\prn{\mc N^{\mc S},\mc S}\ll 1$.

Using \eqref{eq:ProfitLoc}, we define the cost of complexity for system $\cur{\mc N^{\mc S}, \mc S}$ as
\begin{align}\label{eq:CompCost}
 C_\alpha \prn{\mc N^{\mc S},\mc S}=\sum_{\ell \in \mc S}\alpha C_p\prn{\mc N_\ell,\ell} R\prn{\ell}
\end{align}

Additionally, expressions \eqref{eq:ProfitLoc} and \eqref{eq:CompCost} can be easily generalised to the case where the cost of complexity is different for each facility by simply defining independent $\alpha_\ell$ values for each facility $\ell \in \mc S$.

\subsection{The $K$-MedianPlex Problem}

The objective of the decision problem is to find the optimal number and location of facilities, and the associated demand allocation, to maximise the firm’s profit when considering the cost of complexity. 

We first introduce the optimisation problem for a given number $K$ of facilities. We refer to this as the $K$-MedianPlex problem. Using $\mc{S}$ for a given set of open facilities, and $\phi_\ell$ for the fix opening/operating cost of facility $\ell \in \mathcal{S}$, our problem can be written as:

\begin{align}
Z^K_{Plex} = \max_{\mc S \subset \mc N} \quad &   \Pi^\circ \prn{\mc N^{\mc S},\mc S}-\sum_{\ell \in \mc S}\phi_\ell  \label{obj:KPlex}\\
s.t. \quad &  \qquad \abs{\mathcal{S}}=K \label{rest:K}\\
& \qquad \abs{\cur{\ell : i \in \mathcal{N}_\ell}}=1, \quad i \in \mathcal{N} \label{rest:alloc1}
\end{align}
where $\Pi^\circ \prn{\mc N^{\mc S},\mc{S}}$ is given by equation \eqref{eq:ProfitLoc}. Please notice that when using expression $\Pi^{\circ}\prn{\mc N^{\mc S},\mc{S}}$, the local complexity, $C_p\prn{\mc N_\ell,\ell}$, is obtained by using the contribution of each demand node to local demand served by facility $\ell$, i.e. $q^\ell_i=\dfrac{W_i}{\sum_{h\in\mathcal{N}_\ell} W_h}$ for all $i\in\mathcal{N}_\ell$. Constraint \eqref{rest:K}, the cardinal of set $\mc S$, indicates the number of open facilities; and constraints \eqref{rest:alloc1} impose that each demand node should be allocated to only one facility and guarantees that all nodes are allocated.

It can be seen that the $K$-MedianPlex problem is a non-linear variant of a $K$-Median problem with fixed costs:
\begin{align}
 Z^K = \max_{S \subset \mc N} \quad &  \sum_{\ell \in \mathcal S}\bra{ R\prn{\ell} - \phi_\ell } \label{Obj:K-med}\\
s.t. \quad &  \qquad \abs{\mathcal{S}}=K \\
& \qquad \abs{\cur{\ell : i \in \mathcal{N}_\ell}}=1, \quad i \in \mathcal{N} 
\end{align}
 We use $\mc S^K$ to represent the optimal solution to problem \eqref{Obj:K-med} when the number of open facilities is $K$. Solving the $K$-MedianPlex problem for different values of $K$ allows finding the number and location of facilities that maximises the firm's profit under the presence of structural complexity. 

\vspace{12pt}

The $K$-MedianPlex problem presents high computational complexity. It is a non-linear, highly combinatorial problem, hardly solvable for real-sized problems. However, structural complexity is not, in general, a design problem, but a consequence of successive network expansions aimed at reaching higher profitability. This is aligned with the objective of this work which, rather than finding the optimal solution to the $Z^K_{Plex}$ problem, seeks to identify the reasons why network restructuring efforts may fail at reducing operating costs; and to assess the possibility of reducing the burden of complexity in a network. With the second objective in mind, we propose the following $Z^\circ_{Plex}$ function, which can be used to find the net profit of any distribution network $\cur{\mc N^{\mc S},\mc S}$:  

\begin{align}\label{eq:local}
 Z^\circ_{Plex}\prn{\mathcal{N}^{\mc S},\mc S}=\sum_{\ell \in \mc S} \prn{ R\prn{\ell}\prn{1-\alpha C_p\prn{\mc N_\ell, \ell}}-\phi_\ell \str}
\end{align}

Without any loss of generality, in what follows it is assumed that network $\mc S$ has been designed without taking into account the cost of complexity or, alternatively, that it is the product of successive expansions of a smaller network (we refer to this as \textit{organic growth}).

\section{Assessing the impact of complexity}\label{sec:ineffective}

As we have pointed out, one of the advantages of our formulation is to provide practical and intuitive tools to assess the cost of locational complexity, and to propose an explanation to both, the so-called addiction to growth, and the ineffectiveness of restructuring processes suggested by the empirical evidence.

\subsection{The cost of complexity}

It has been said that complexity introduces hidden costs that cripple the capacity of the firm for generating profits \citep{WuFriEfs2007, BozWarFly2009}. In location decisions, additional costs required for connecting, supporting, and supplying the distribution network, may outweigh the increased revenues generated from the access to new markets and customers.

\begin{figure}[!ht]
	\centering
	\mbox{\includegraphics[trim={0 0 0 0},width=8cm]{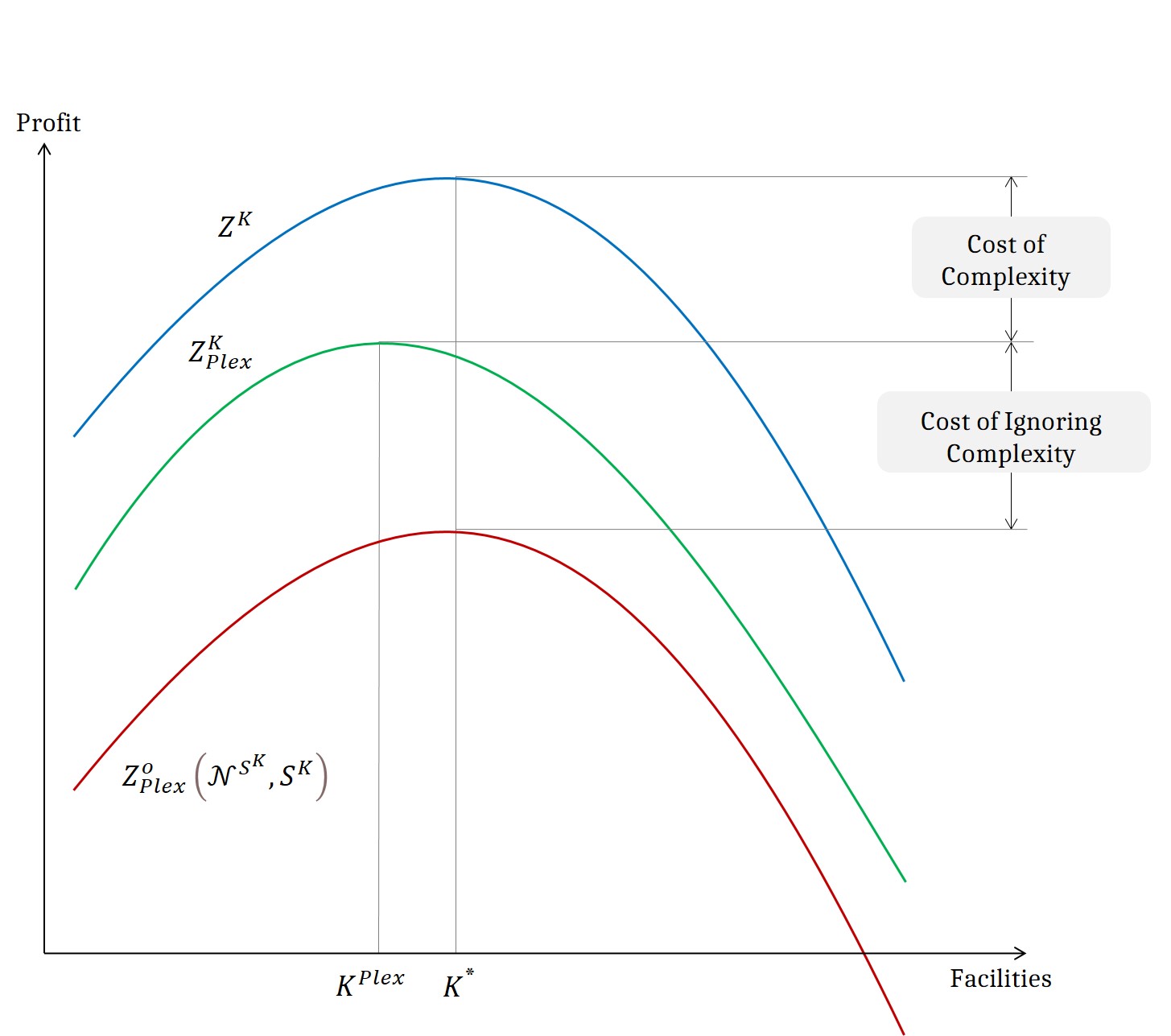}}
	\caption{Cost of Complexity }\label{fig:compcost}
\end{figure}

The $K$-MedianPlex problem formulation, represented by equations \eqref{obj:KPlex} to \eqref{rest:alloc1}, together with equation  \eqref{eq:local}, allow us to provide a visual characterisation of the burden of locational complexity. The top line in Figure \ref{fig:compcost} represents the forecasted profits implied from standard facility location models, in particular the $K$-Median problem, for different numbers of facilities. The bottom curve represents the actual profits. This curve is obtained by evaluating the current network configuration in equation \eqref{eq:ProfitLoc}. Finally, the intermediate curve represents the optimal solution to the $K$-MedianPlex problem. The distances between these curves represent two different components of the cost imposed by locational complexity on the company. First, there is an unavoidable cost of complexity. It is the cost imposed by the intricate interrelations existing in any distribution network. This cost is represented by the distance between the upper and the intermediate curves. The second one is the cost of ignoring the existence of complexity during network design, and it is represented by the distance between the intermediate and the bottom curves.  This discussion will be recalled in Section \ref{sec:algor} when proposing strategies aimed at alleviating the burden of complexity in network design.

\subsection{The addiction to growth}

Regarding the addiction to growth, complexity-driven costs are usually hidden in early design stages. This implies that the observed profit may be, in practice, smaller than the one forecasted during the planning process. If revenues are as expected, this may give the management the wrong impression that the optimal number of facilities may indeed be larger than initially planned. The natural reaction might be to further increase the number of facilities, aiming at reaching the profit levels forecasted by the analytical models. If this actually leads to an increase in profit, but the observed figures remain below the forecasted ones, the same reasoning can lead the company to unnecessary expansions in the network in the quest for the higher profit levels promised by the analytical models.  

\begin{figure}[!ht]
	\centering
	\includegraphics[trim={0 2cm 0cm -0.5cm}, height=6cm]{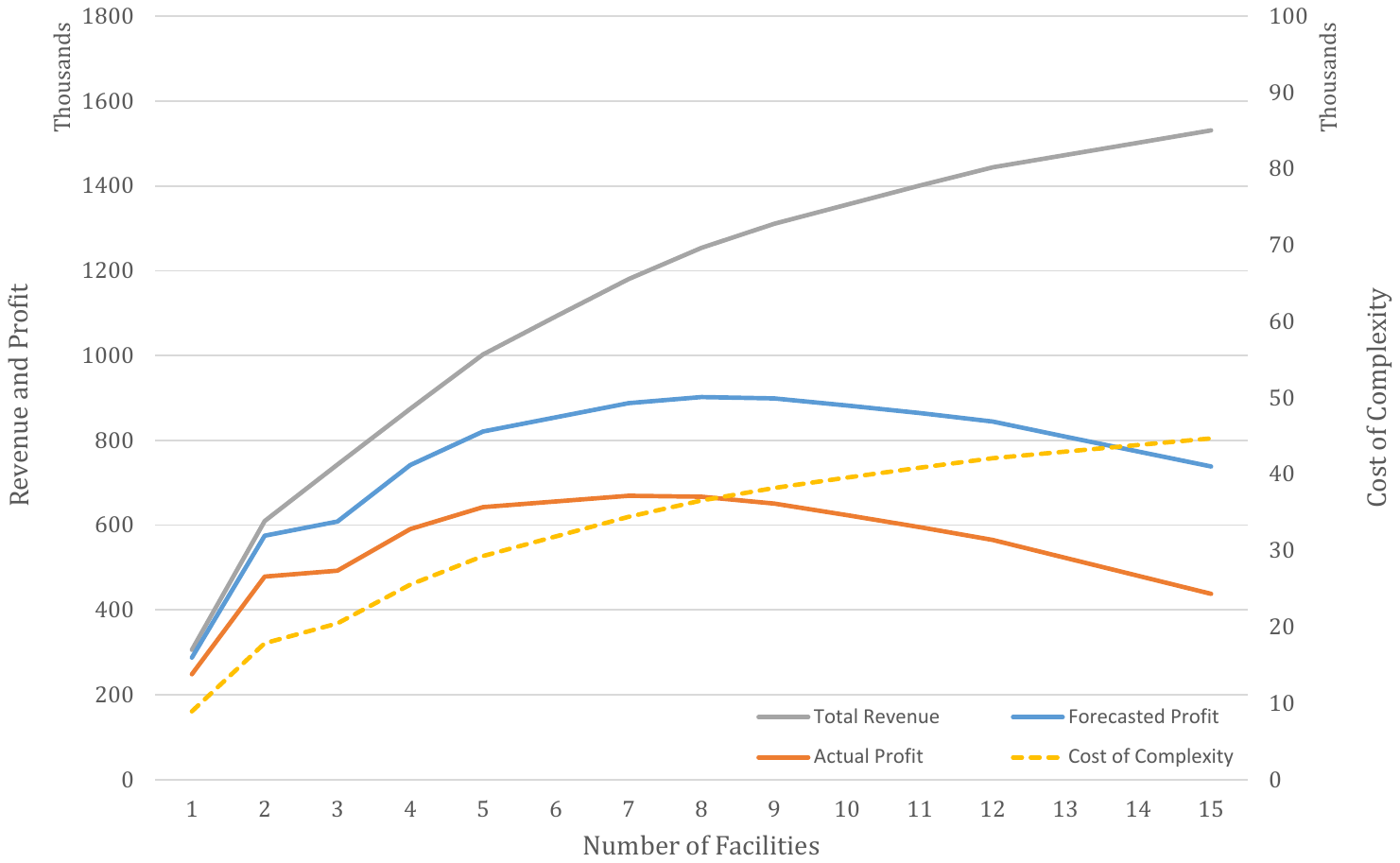}
	\caption{\cm{Addiction to growth in network design.}}\label{CostComp}
\end{figure}

A simple numerical exercise, using the $K$-Median Plex objective function \eqref{obj:KPlex} with $\Pi^\circ\prn{\mc N^{\mc S},\mc S}$ given by equation \eqref{eq:ProfitLoc}, and approximating the complexity cost $\alpha$ to 2.5\% of the profits as suggested by \citet{MenRui2020}, provides support to this claim. Figure \ref{CostComp} depicts the solution to a number of $K$-Median problems, ranging from 1 to 15 facilities. The upper curve represents the total revenue associated to each solution. The forecasted and actual profits are represented by the two concave curves (left-hand side axis). The cost of complexity, the distance between the two profit curves, is represented by the dashed curve (right-hand side axis). It can be seen that a decision maker that focuses solely on revenue, may have a strong incentive to keep expanding his network in the quest for higher profit levels, ignoring that the company has grown beyond its optimal size.

\subsection{Restructuring and complexity} \label{sec:FutRes}

Data from different retailers and banks evince that, despite large network reduction and restructuring efforts, their operating costs have either experienced only a minor reduction or remained almost constant. Straightforward manipulation of the \textit{pars}-Complexity measure may provide insight for understanding this phenomenon. 

A common feature of all these companies is that their service is universal, i.e. it targets all the regional demand and no client can be excluded. Therefore, irrespectively of the number of facilities, the target market remains the same. The reader can think, for example, of a convenience-stores chain that has three shops in a neighbourhood. In the absecne of competitors, closing one of their facilities will trigger a re-distribution of clients across the facilities that remain open, but the demand will remain fairly the same. The closing of bank branches, where the clients are simply allocated to a different branch is another example of this.

In order to bring this into the locational complexity framework, consider a centrally managed system $\mc L^1$ consisting of $K$ independent subsystems or regional distribution centres $\mc S^1_1, \ldots, \mc S^1_K$. Each of these centres contributes a proportion $p^1_\ell, \, \ell= 1,\ldots, K$ of the total sales. These contributions are collected in vector $\mb p^1=\prn{p^1_1,\ldots,p^1_K}$. Likewise, each distribution centre $\mc S^1_\ell$ has associated a number $H_\ell$ of retailers, each contributing a proportion $q^{1,\ell}_h$ of the centre's total sales. The weights vector of centre $\ell$ is represented by $\mb q^{1,\ell}=\prn{q^{1,\ell}_1,\ldots,q^{1,\ell}_{H_\ell}}$. The \textit{pars}-Complexity of this system, represented by $C^1_p=C_p\prn{\mc L^1}$, can be readily computed using the right hand side expression in \eqref{eqn:CompTot}. Consider an alternative configuration of the network, system $\mc L ^2$, consisting of $M$ independent regional distribution centres, each of them serving a number $ G_m$ of retailers. The corresponding weights vectors are $\mb p^2$ and $\mb q^{2,m}, m=1,\ldots,M$. The \textit{pars}-Complexity is represented by $C^2_p=C_p\prn{\mc L^2}$. Assume also that $M \ll K$, and $N=\sum_{\ell=1}^K H_\ell = \sum_{m=1}^M G_m$ is the total number of retail centres. 

The following proposition uses the convention $\oplus_{i=1}^n \mb g_i= \mb g^1 || \mb g^2 ||\cdots || \mb g^n$ for representing the concatenation of vectors $\mb g^1$ to $\mb g^n$. 

\begin{proposition} \label{Prop:1}
 Consider systems $\mc L^1$ and $\mc L^2$. If $\sum_{\ell=1}^L H_\ell = \sum_{m=1}^M G_m$ and $\oplus_{\ell=1}^K \prn{\mb p^1_{\ell}\mb q^{1,\ell}}= \oplus_{m=1}^M \prn{\mb p^2_{\prn{m}}\mb q^{2,m}}$, i.e. the contribution of each retailer to the total demand remains the same, then the \textit{pars}-Complexity values for both networks is the same, i.e. $C^1_p\prn{\mb g^1} = C^2_p\prn{\mb g^2}$.
 
\noindent \begin{proof}
  Let $\mb g^1=\oplus_{\ell=1}^K \prn{\mb p^1_{\ell}\mb q^{1,\ell}}$ and $\mb g^2=\oplus_{\ell=m}^M \prn{\mb p^2_{m}\mb q^{2,m}}$. Using equation \eqref{eq:Cp} and given the fact that $\mb g^1 = \mb g^2$ it comes straightforwardly that
  \[C_p \prn{\mb g^1}=\sum_{i=1}^{N} g^1_i =\sum_{i=1}^{N} g^2_i =  C_p \prn{\mb g^2}   \]
  i.e. the structural complexity value of both systems is identical.
 \end{proof}
\end{proposition}

Proposition \ref{Prop:1} provides an intuitive result that finds its foundation in the Second law of Thermodynamics (roughly speaking, the total entropy in a closed system can never decrease): if structural complexity is understood --as \citet{RuiMenAmr2019} suggest- as an entropic process, then the level of complexity faced by a company will be the same as long as the target market remains unchanged, irrespective of the actual network configuration.

Therefore, as long as their market remains the same, companies that provide in-store services will hardly attain any reduction in complexity-related costs by simply eliminating intermediate distribution centres or local facilities. Moreover, the elimination of one centre (and the corresponding absorption of its associated demand by one or more facilities) may increase the total complexity of the facilities that remain open. 

These results also suggest that the only mechanism that may bring a reduction in \textit{pars}-Complexity, and therefore in operating costs, is a reduction in its market share. For the interested reader, \citep{RuiMenAmr2019} provides a mechanism for computing the market contribution at which the elimination of a certain demand node will surely reduce supply chain's structural complexity.

In the following section we present three alternative mechanisms (two of them involving demand reduction) for reducing the burden of locational complexity  in the supply chain.

\section{Addressing locational complexity}\label{sec:algor}

In this section we adopt the view of a firm that is already suffering the burden of complexity and seeks to identify mechanisms for reducing the cost of structural complexity in their network. Three alternative approaches are considered: in the first case, \textit{network rebalancing}, the firm aims at reducing the cost of complexity at the local level by means of a reassignment of the markets across different facilities. In the second one, the \textit{network rationalisation} approach, the firm abandons non-profitable complexity-bearing markets. Finally, in the third case, \textit{network reduction}, we assume that the firm closes a number of facilities and has the option to either completely abandon the associated market(s), or to reallocate them to a different facility that remains open.

The efficiency of the proposed approaches is assessed in Section \ref{sec:NumEx} by evaluating different network configurations, characterised by a set of facilities $S$ and their corresponding allocation sets $\mathcal{N}^{\mc S}$, on each of the equations $Z^K$, and $Z^\circ_{Plex}$.

\begin{figure}[!h]
	\centering
	\mbox{\includegraphics[width=9cm]{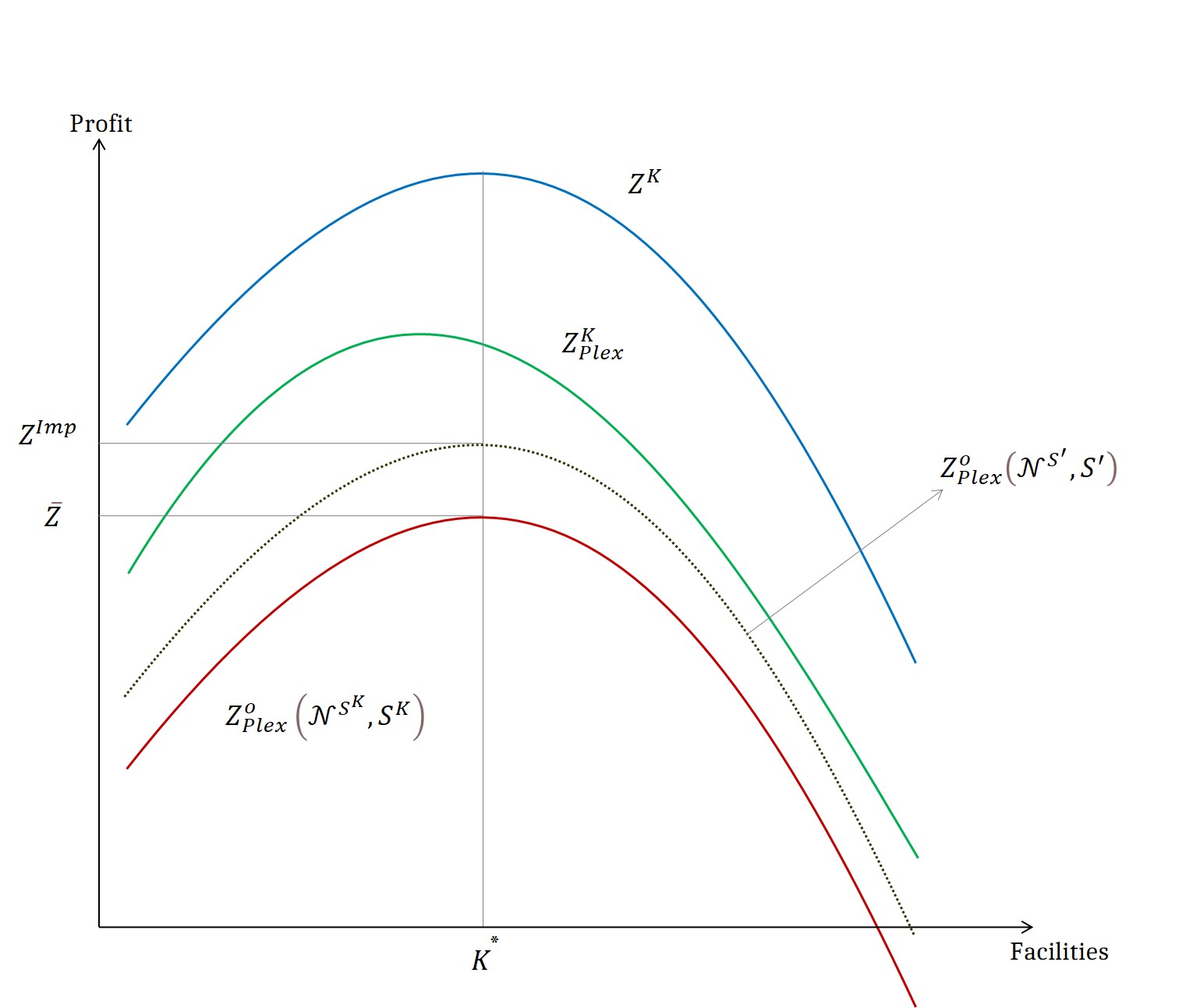}}
	\caption{\cm{Forecasted, actual and improved profits for different numbers of facilities.} }\label{CostComp2}
\end{figure}

This rationale is illustrated in Figure \ref{CostComp2}. The top and bottom curves represent, respectively, the profit $Z^K$ implied by the solution $\cur{\mc N^{{\mc S^K}},\mc S^K}$ to the $K$-Median problem for different values of $K$; and the actual profit observed of the company, $Z_{Plex}^\circ$, obtained by evaluating the different $K$-Median solutions in equation \eqref{eq:local}. Likewise, the $Z^K_{Plex}$ curve represents optimal solution to the corresponding $K$-MedianPlex problem for each value of $K$. The $Z^\circ_{Plex}(\mathcal{N}^{\mc S'}, \mc S')$ curve, represents the improvement that could -at least hypothetically- be attained by means of any of the three strategies proposed in this section. This curve is obtained by evaluating alternative solutions $\cur{\mc N^{\mc S'},\mc S'}$ on equation \eqref{eq:local}. Finally, $K^*$ represents the profit maximising number of facilities in the $K$-Median formulation (this number does not necessarily coincide with the optimal number of facilities prescribed by the $Z^K_{Plex}$ formulation); $\ov Z$ provides the actual profit associated to the $\mc S^{K^*}$ solution, and $Z^{Imp}$ indicates the profit attained by the improvement approaches for a number $K^*$ of facilities. The objective of the approaches proposed in this section, is finding a network configuration such that the $Z^\circ_{Plex}(\mathcal{N}^{\mc S'},\mc S')$ curve gets as close as possible to $Z^K_{Plex}$ for any given number of facilities.

Sections \ref{sec:rede} to \ref{sec:NetRed} describe in detail each of the proposed approaches and provide algorithmic procedures for identifying potential improvements on the network.

\subsection{Network rebalancing} \label{sec:rede}

Standard facility location and network design models ignore the hidden costs of channel and market proliferation. This strategy results in oversized networks whose profitability is often overestimated. This situation frequently leads companies to look for strategies aimed at improving their network's performance \citep{FisGauKle2017}.  In order to address this problem, in particular in situations where local facilities still appear profitable, it may seem intuitive to reduce the burden carried by the facilities that bear larger demand, and to cede market share to other facilities which currently carry a lighter load. We refer to this type of strategy as \textit{network rebalancing}.

To assess the effectiveness of such policies, in this section we propose a strategy that, starting from a solution to the $K$-Median problem characterised by the pair $\left\{\mathcal{N}^{\mc S^K}, \mc S^K\right\}$, reassigns demand nodes to different facilities exploiting trade-offs between increased transportation costs and reduced complexity penalties. The reassignment is conducted by sequentially selecting nodes in decreasing order of the total shipment cost to the second closest facility. If reallocating a node to its second closest facility brings an improvement in the objective function, the node is reallocated; otherwise, it remains assigned to its original facility. After demand reallocation, the algorithm takes one further improvement step by solving a $1$-Median problem for each of the available facilities and their associated nodes.

Pseudo-code for this reallocation and re-centring procedure is presented in \ref{alg:ZPlex}. A numerical assessment of the efficiency of this strategy is presented in Section \ref{sec:NumEx}.

\subsection{Network rationalisation} \label{sec:ratio}

Locational complexity is not typically result of network design, but a problem that arises from organic growth, i.e. from successive network expansions aimed at achieving growth in profits. Thus, in order to reduce complexity, firms may consider abandoning certain markets and concentrating their services in the most profitable regions (although this approach often finds resistance under the rationale that lost sales will affect profit negatively). To address this issue, we propose a \textit{network rationalisation} strategy consisting of sequentially removing those demand nodes for which the savings from reducing complexity and transportation costs, offset the lost revenue and lead to an increase in profits. 

The pseudo-code is presented in \ref{alg:DemRed}. The algorithm starts from a given network configuration $\left\{\mathcal{N}^{\mc S^K}, \mc S^K\right\}$ and takes, for each given facility $\ell \in \mc S^K$, a set with an arbitrary number of the most distant demand nodes in $\mathcal{N}_\ell$. We refer to this set as $\Theta_\ell$. For each node $i \in \Theta_\ell$, the algorithm computes the value of the $Z_{Plex}^\circ$ function (given by equation \eqref{eq:local}), represented by $Z_{Plex}^\circ \left(\mathcal{N}^{\mc S^K}-\left\{i\right\},\mc S^K \right)$,  which results from eliminating $i$ from $\mathcal{N}_\ell$ and computes the ratio
\[\lambda(i)=\frac{Z_{Plex}^\circ\left(\mathcal{N}^{\mc S^K}-\left\{i\right\},\mc S^K \right)}{\widetilde{Z}_{Plex}}-1 \, , \]
where $\widetilde{Z}_{Plex}$ is the best available value of the objective function. This value is initialised as $ \rule{0pt}{12pt} \widetilde{Z}_{Plex}=Z_{Plex}^\circ\left(\mathcal{N}^{\mc S^K},\mc S^K \right)$. The node that returns the largest (positive) value of $\lambda\left(i\right)$ is eliminated, and the procedure repeated until no further improvement in the objective function can be attained from facility $j$. This strategy is numerically assessed in Section \ref{sec:NumEx}.

\subsection{Network reduction}\label{sec:NetRed}

In Section \ref{sec:FutRes} we argued that only when the elimination of facilities is accompanied by a reduction on market share, the company can aim at reducing  complexity-related costs. This situation is addressed by a simple procedure consisting of sequentially eliminating facilities (and their associated demand nodes) as long as an improvement in the network's gross profit can be attained, i.e. facility $\ell \in \mc S' \subseteq \mc S^K$ will be eliminated if and only if $Z_{Plex}^\circ \prn{\mathcal{N'}^{\mc S'-\{ \ell \}}-\mathcal{N}_\ell,\mc S'-\{ \ell \}} > Z_{Plex}^\circ \prn{\mathcal{N}^{\mc S'},\mc S'}$, where $\mathcal{N'} \subseteq \mathcal{N}^{\mc S^K}$, and  $\cur{\mathcal{N}^{\mc S ^ K},\mc S^K}$ is the initial network.  On each iteration of the algorithm, we eliminate the facility whose removal brings the largest improvement in the objective function. The routine, whose pseudo-code is provided in \ref{alg:NetRed}, stops when no further improvement can be obtained by eliminating one more facility.

For of completeness, we also assess the possibility of only eliminating the facility, reallocating the demand nodes among the facilities that remain open. In such case, the gross profit is given by $Z_{Plex}^\circ \prn{\mathcal{N}^{\mc S'-\{ \ell \}},\mc S'-\{ \ell \}}$, where $\ell \in \mc S'$ represents the eliminated facility. This strategy is evaluated in Section \ref{sec:NumEx}.

\section{Numerical experiments}\label{sec:NumEx}

The objective of this section is to assess the effectiveness of alternative strategies for reducing the burden of structural complexity in an oversized network.
The experimental setting is a distribution network including cities with more than 50 thousand inhabitants in Spain.

The heuristics were tested on networks designed over the 125 candidate cities. Nearly 10 thousand examples were conducted. Each instance consists of a sequence of $K$-Median problems, with the number of facilities, $K$, ranging from $2$ to $9$. For each network, we recorded the corresponding solution $\cur{\mc N^{\mc S^K},\mc S^K}$; the associated values of the $K$-Median objective function, $Z^K$; and the observed profit, $Z^\circ_{Plex} \prn{\mc N^{\mc S^K},\mc S^K}$. Function $Z^\circ_{Plex}$  was evaluated for different values of parameter $\alpha$ taken from the set $\{0.025,0.05, \allowbreak 0.075, \allowbreak 0.1,0.125,0.15\}$, and for values of $\gamma$ in set  $\{8.3, 16.6, 33.3, 66.6,100,200,\allowbreak 400\}$, representing cents/km per ton. Additionally, three different values for the facilities operating costs, $\phi$, were considered, $\cur{50k, 70k, 80k}$.

In an alternative specification, we assumed that unit transport costs from the main distribution centre to the regional facilities, $\varrho$, may differ from the ones incurred by the regional facilities when serving their nodes, $\gamma$. In such case, equation \eqref{eq:revenue} becomes
\begin{align}\label{eq:revalt}
R\left(\ell \right)=\sum_{i\in\mathcal{N}_\ell}{\left(r-\varrho d\prn{c,\ell}-  \gamma d\prn{\ell,i}\right)\, W_i}, \quad & \ell \in \mathcal{S} 
\end{align}
where $d\prn{c,\ell}$ represents the distance from the main distribution centre to regional facility $\ell$, and $\varrho$ represents the transport costs per ton and km in the first leg. Notice that for $\varrho=0$, equation \eqref{eq:revalt} reduces to \eqref{eq:revenue}. In our numerical experiments, parameter $\varrho$ is taken from the set $\cur{8.3, 16.6,33.3,66.6,100,200}$. Additionally, we assume that $\varrho < \gamma$. \jp{All distances are by road and were obtained using Google Maps API. Demand values were obtained based on information for year 2017 provided by the Spanish National Institute of Statistics (INE), and scaled using the expression: $W_i=500 \ln Population_i $.}

All instances were ran in Matlab \textregistered\ using the IBM ILOG CPLEX \textregistered\ connector. In all cases the execution time was below 15 secs.  Given the strategic nature of the problem, and the fact that the main purpose of the heuristics is to assess the efficiency of alternative solutions that could be taken in practice for addressing locational complexity problems, we consider that the heuristics' computational performance (in terms of both, time, and memory usage) is of little to no relevance for this manuscript and, therefore, no further details are provided.

Table \ref{tab:kvalues} presents the profit-maximising number of facilities, $K^*$ and $K^\circ_\alpha$, under objectives $Z^K$ and $Z^\circ_{Plex}$ respectively, for different values of parameters $\alpha, \gamma, \phi$ and $\rho$. The revenue associated to different values of $\alpha$ was computed using equation \eqref{eq:revalt}.

The three different approaches described in Section \ref{sec:algor} were deployed for restructuring the network associated to each $K$-Median instance. Our results suggest that whereas some improvements in profit can be attained by means of the \textit{network rebalancing} strategy, the impact of \textit{network rationalisation} on profits is minimal. It is also observed that significant improvements in gross profit can be attained by means of \textit{network reduction}, in particular when the demand associated with an eliminated node is reasigned among those facilities that remain open.  As it may be expected, the improvement routines return higher improvements for large values of the complexity cost parameter $\alpha$. A detailed discussion of our results is presented below.

\begin{sidewaystable}
\centering
\scriptsize
\begin{tabular}{c|c|cccccc|cccccc|cccccc} \hline
 
 \multicolumn{2}{c|}{\st }				& \multicolumn{6}{c|}{$\phi=50000$}											&	\multicolumn{6}{c|}{$\phi=70000$}											&	\multicolumn{6}{c}{$\phi=80000$}											\\	\cline{3-20}
\multicolumn{2}{c|}{\st} & \multicolumn{18}{c}{$\gamma$}          \\ \hline
\st	$\rho$&		&	400	&	200	&	66.6	&	33.3	&	16.6	&	8.3	&	400	&	200	&	66.6	&	33.3	&	16.6	&	8.3	&	400	&	200	&	66.6	&	33.3	&	16.6	&	8.3	\\	\hline
\multirow{4}{*}{\rotatebox{90}{200}}	&	\st $K^*$	&	3	&		&		&		&		&		&	3	&		&		&		&		&		&	3	&		&		&		&		&		\\	
\st		&	$K^\circ_{0.075}$	&	5	&		&		&		&		&		&	3	&		&		&		&		&		&	3	&		&		&		&		&		\\	
\st		&	$K^\circ_{0.1}$	&	5	&		&		&		&		&		&	5	&		&		&		&		&		&	3	&		&		&		&		&		\\	
\st		&	$K^\circ_{0.125}$	&	5	&		&		&		&		&		&	5	&		&		&		&		&		&	5	&		&		&		&		&		\\	\hline

\multirow{4}{*}{\rotatebox{90}{66.6}}	&\st	$K^*$	&	6	&	3	&		&		&		&		&	5	&	3	&		&		&		&		&	4	&	2	&		&		&		&		\\	
\st		&	$K^\circ_{0.075}$	&	7	&	6	&		&		&		&		&	6	&	4	&		&		&		&		&	5	&	4	&		&		&		&		\\	
\st		&	$K^\circ_{0.1}$	&	7	&	6	&		&		&		&		&	6	&	5	&		&		&		&		&	5	&	4	&		&		&		&		\\	
\st		&	$K^\circ_{0.125}$	&	7	&	7	&		&		&		&		&	6	&	6	&		&		&		&		&	6	&	5	&		&		&		&		\\	\hline

\multirow{4}{*}{\rotatebox{90}{33.3}}& \st	$K^*$	&	6	&	3	&	2	&		&		&		&	5	&	3	&	2	&		&		&		&	5	&	3	&	2	&		&		&		\\	
\st		&	$K^\circ_{0.075}$	&	7	&	6	&	5	&		&		&		&	6	&	4	&	5	&		&		&		&	5	&	4	&	3	&		&		&		\\	
\st		&	$K^\circ_{0.1}$	&	7	&	7	&	5	&		&		&		&	6	&	6	&	5	&		&		&		&	5	&	4	&	5	&		&		&		\\	
\st		&	$K^\circ_{0.125}$	&	7	&	7	&	5	&		&		&		&	6	&	6	&	5	&		&		&		&	6	&	6	&	5	&		&		&		\\	\hline

\multirow{4}{*}{\rotatebox{90}{16.6}}	&	\st$K^*$	&	6	&	4	&	2	&	2	&		&		&	5	&	3	&	2	&	2	&		&		&	5	&	3	&	2	&	2	&		&		\\	
\st		&	$K^\circ_{0.075}$	&	7	&	6	&	6	&	5	&		&		&	6	&	5	&	4	&	5	&		&		&	5	&	5	&	4	&	3	&		&		\\	
\st		&	$K^\circ_{0.1}$	&	7	&	7	&	6	&	5	&		&		&	6	&	5	&	5	&	5	&		&		&	5	&	5	&	4	&	5	&		&		\\	
\st		&	$K^\circ_{0.125}$	&	7	&	7	&	7	&	5	&		&		&	7	&	6	&	6	&	5	&		&		&	5	&	5	&	6	&	5	&		&		\\	\hline

\multirow{4}{*}{\rotatebox{90}{8.3}}	&	\st $K^*$	&	6	&	5	&	2	&	2	&	2	&		&	5	&	3	&	2	&	2	&	2	&		&	5	&	3	&	2	&	2	&	2	&		\\	
\st		&	$K^\circ_{0.075}$	&	7	&	6	&	5	&	6	&	5	&		&	6	&	5	&	5	&	4	&	5	&		&	5	&	5	&	4	&	4	&	3	&		\\	
\st		&	$K^\circ_{0.1}$	&	7	&	7	&	7	&	6	&	5	&		&	6	&	5	&	5	&	5	&	5	&		&	5	&	5	&	5	&	4	&	5	&		\\	
\st		&	$K^\circ_{0.125}$	&	7	&	7	&	7	&	7	&	5	&		&	7	&	6	&	5	&	6	&	5	&		&	5	&	5	&	5	&	6	&	5	&		\\	\hline
																																									
\multirow{4}{*}{\rotatebox{90}{\parbox{1.3cm}{\centering W/o Dist. Costs}}}	&	\st$K^*$	&	6	&	5	&	2	&	2	&	2	&	2	&	5	&	3	&	2	&	2	&	2	&	2	&	5	&	3	&	2	&	2	&	2	&	2	\\	
\st		&	$K^\circ_{0.075}$	&	7	&	7	&	5	&	5	&	5	&	5	&	6	&	5	&	5	&	4	&	4	&	3	&	5	&	5	&	3	&	3	&	3	&	3	\\	
\st		&	$K^\circ_{0.1}$	&	7	&	7	&	7	&	7	&	7	&	7	&	6	&	5	&	5	&	5	&	5	&	5	&	5	&	5	&	5	&	5	&	5	&	5	\\	
\st		&	$K^\circ_{0.125}$	&	7	&	7	&	7	&	7	&	7	&	7	&	7	&	6	&	5	&	5	&	5	&	5	&	6	&	5	&	5	&	5	&	5	&	5	\\	\hline

\end{tabular}
 \caption{Profit maximising number of facilities for the $K$-Median and the associated $Z^\circ_{Plex}$ formulation} \label{tab:kvalues}
 \end{sidewaystable}

\subsection{Network rebalancing}

This strategy attempts to reduce the cost of complexity by balancing the complexity value of the different facilities through a reallocation of nodes. After reallocating nodes, the algorithm seeks an improvement of the network by recentering some facilities within their own network. 

Table \ref{fig:Redis} shows the results of 270 experiments. The table consists of three blocks associated with different values of the operation costs parameter ($\phi$). Each of these blocks is divided into three additional blocks, each of them associated with a value of the complexity cost parameter ($\alpha$). Finally, each of the nine sub-blocks provides the results obtained for networks sizes from 5 to 9 facilities. All the cases was solved for six alternative values of the transport cost parameter (columns). The values in the table represent the relative percentage improvement ($\Delta Z^\circ_{Plex}$) attained by the routine with respect to the $Z^\circ_{Plex}\prn{\mathcal N^{\mc S},\mc S }$ value given by \eqref{eq:local}.

The results suggest that the \textit{network rebalancing} heuristic performs better for lower transport costs, higher operation costs, higher values of the complexity parameter, and for networks with a larger number of facilities. This confirms that \textit{network rebalancing} efforts may be more effective for highly over-dimensioned networks. Our results also suggest that the larger the transport costs, the smaller the potential impact on profit of this strategy.

\begin{figure}
\hspace{-12pt}\begin{subfigure}{0.6\textwidth}
\includegraphics[width=6.5cm]{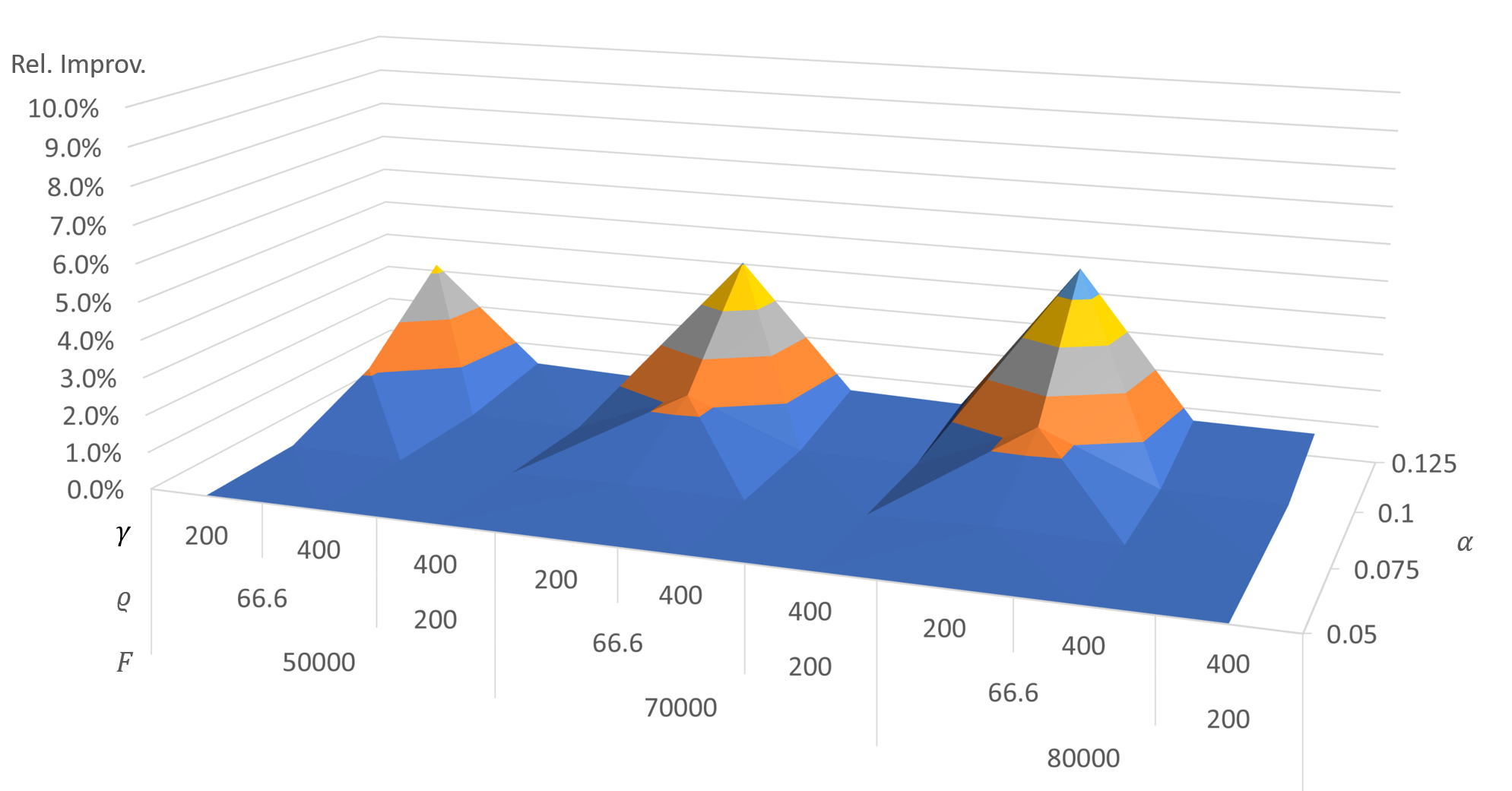} 
\caption{K=5}
\end{subfigure}
\begin{subfigure}{0.6\textwidth}
\includegraphics[width=6.5cm]{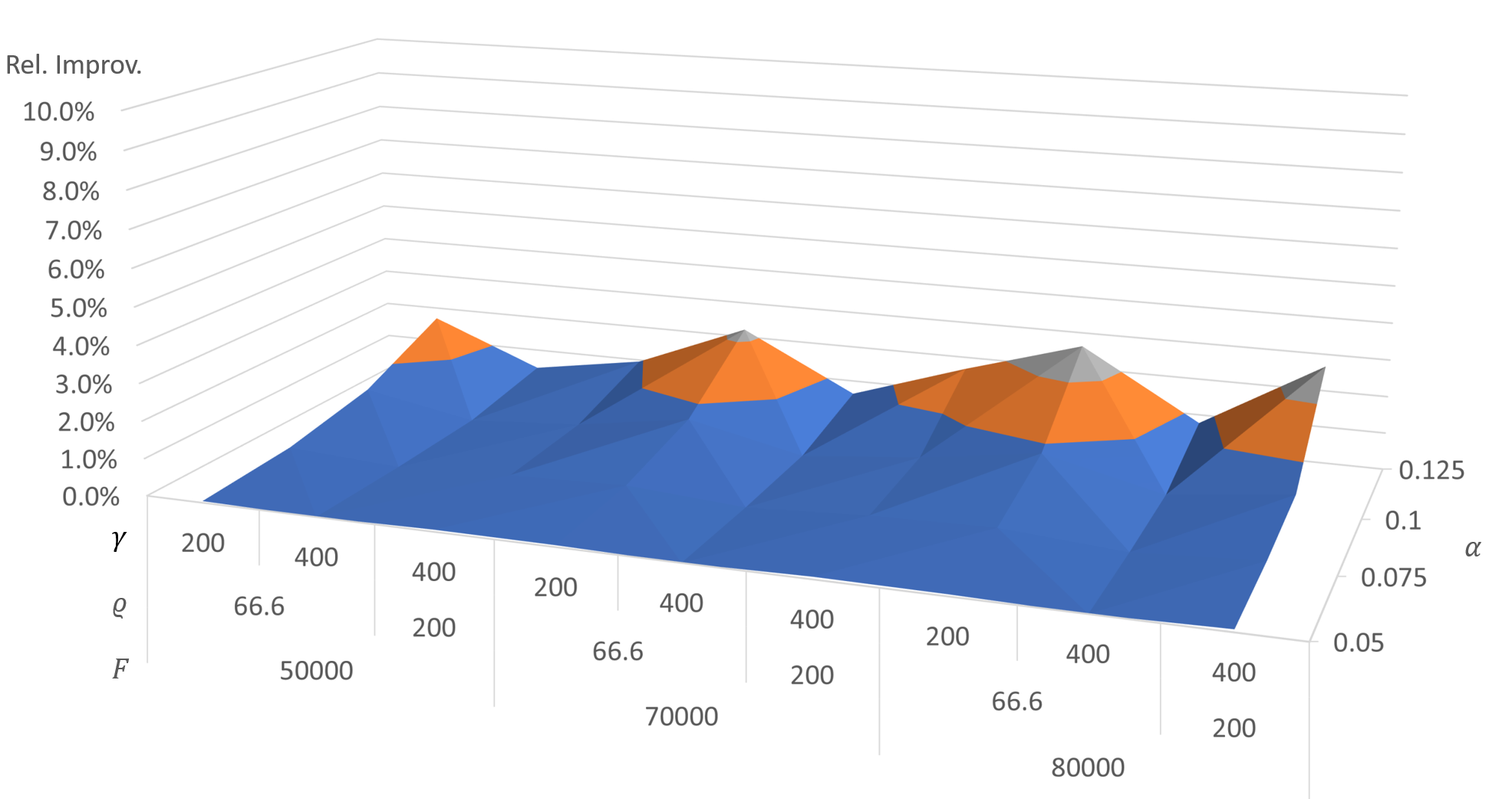} 
\caption{K=6}
\end{subfigure} \\

 \hspace{-12pt}\begin{subfigure}{0.6\textwidth}
\includegraphics[width=6.5cm]{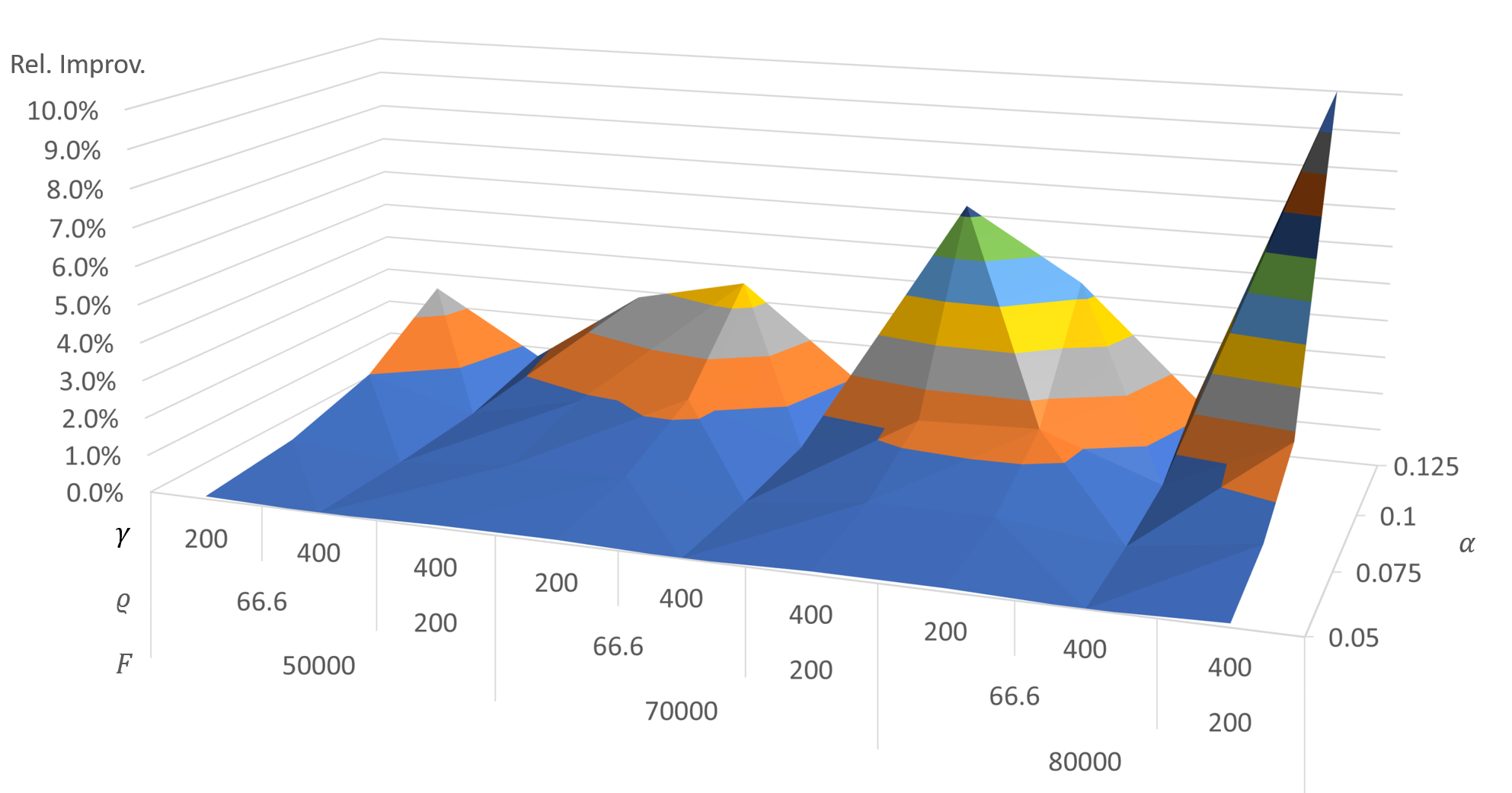} 
\caption{K=7}
\end{subfigure}
\begin{subfigure}{0.6\textwidth}
\includegraphics[width=6.5cm]{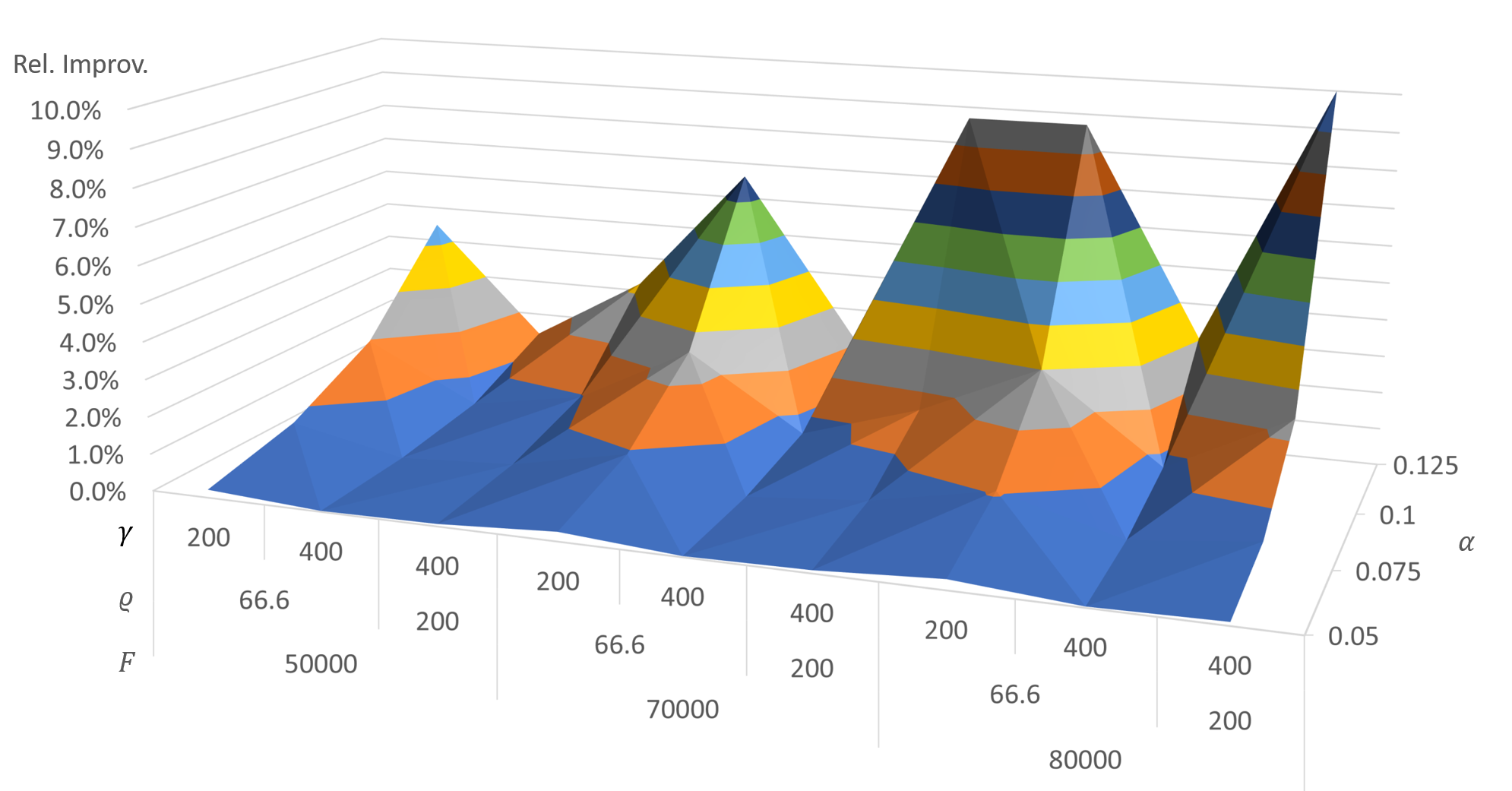} 
\caption{K=8}
\end{subfigure}
\caption{Relative improvement with the rebalancing strategy for different values of K}\label{fig:RedStrPlots}
\end{figure}

The results for the case involving distribution costs, with revenue function represented by equation \eqref{eq:revalt}, are presented in Table \ref{fig:RedisDist} and in Figure \ref{fig:RedStrPlots}. The results for value $\varrho=8.3$ of the distribution cost, presented in Table \ref{fig:RedisDist} are consistent with what is observed in Table \ref{fig:Redis}. However, it is hard to identify a clear pattern for different values of the distribution cost. This can be better appreciated in Figure \ref{fig:RedStrPlots}, where we present the results of a number of rebalancing experiments, associated with different values of the four modelling parameters ($\alpha, \gamma, \varrho$ and $\phi$). The plots suggest that the rebalancing strategy is more efficient for larger values of K and higher operation costs.

It is important to remember that, given that the number of facilities and their location remains constant, all improvements are due to savings in distribution costs attained by reallocating nodes. Moreover, the total system's complexity remains constant. The impact of the rebalancing strategy on the cost of complexity is limited to marginal changes due to the reallocation of demand nodes, which rebalances the distribution of complexity among facilities. This is shown in Table \ref{fig:RedisComp}. The table shows the change in complexity costs ($\Delta C_\alpha $) obtained in different instances of the numerical experience. Given that the results are the same irrespective of the facility operation costs ($\phi$), we omit this information in the table. It can be observed that in all cases, this strategy attains reductions in the cost of complexity. Moreover, the savings are steeper when the initial network is significantly over dimensioned (9 initial facilities in this case). Table \ref{fig:RedisDistComp} depicts the results obtained by this strategy when distribution costs of $\varrho = 8.30$ are introduced. As it happens in Table \ref{fig:RedisComp}, the results are invariant to changes in the facility operation costs. The results obtained for other values of $\varrho$ follow the same pattern and are therefore omitted.

\begin{table}
\scriptsize
\centering
 \begin{tabular}{c|c|c|cccccc} \hline
	\multicolumn{3}{c|}{Parameters}	& \multicolumn{6}{c}{$\gamma$}	\\	\hline
\st $\phi$ &	$\alpha$	& K &	400.0	&	200.0	&	66.7	&	33.3	&	16.7	&	8.3	\\	\hline
\st 	\multirow{15}{*}{\rotatebox{90}{50000}} & \multirow{5}{*}{\rotatebox{90}{0.075}}	
&	5	&	0.168\%	&	0.193\%	&	0.363\%	&	0.439\%	&	0.475\%	&	0.511\%	\\	
\st	&	&	6	&	0.055\%	&	0.164\%	&	0.397\%	&	0.555\%	&	0.697\%	&	0.780\%	\\	
\st	&	&	7	&	0.040\%	&	0.232\%	&	0.434\%	&	0.545\%	&	0.657\%	&	0.736\%	\\	
\st	&	&	8	&	0.093\%	&	0.212\%	&	0.657\%	&	0.837\%	&	0.923\%	&	1.002\%	\\	
\st	&	&	9	&	0.253\%	&	0.479\%	&	1.203\%	&	1.474\%	&	1.698\%	&	1.817\%	\\ \cline{2-9}

\st &	\multirow{5}{*}{\rotatebox{90}{0.1}}	
&	5	&	0.354\%	&	0.479\%	&	0.727\%	&	0.790\%	&	0.862\%	&	0.907\%	\\	
\st	&	&	6	&	0.215\%	&	0.431\%	&	0.856\%	&	1.109\%	&	1.280\%	&	1.357\%	\\	
\st	&	&	7	&	0.248\%	&	0.513\%	&	0.846\%	&	1.063\%	&	1.193\%	&	1.294\%	\\	
\st	&	&	8	&	0.221\%	&	0.564\%	&	1.274\%	&	1.468\%	&	1.623\%	&	1.733\%	\\	
\st	&	&	9	&	0.537\%	&	1.243\%	&	2.359\%	&	2.723\%	&	2.957\%	&	3.107\%	\\	
	\cline{2-9}

	\st &	\multirow{5}{*}{\rotatebox{90}{0.125}}	
&	5	&	0.938\%	&	1.065\%	&	1.334\%	&	1.480\%	&	1.566\%	&	1.576\%	\\	
\st	&	&	6	&	0.468\%	&	0.995\%	&	1.722\%	&	2.106\%	&	2.316\%	&	2.460\%	\\	
\st	&	&	7	&	0.657\%	&	1.079\%	&	1.620\%	&	1.920\%	&	2.111\%	&	2.236\%	\\	
\st	&	&	8	&	0.610\%	&	1.472\%	&	2.427\%	&	2.641\%	&	2.840\%	&	2.945\%	\\	
\st	&	&	9	&	1.340\%	&	3.034\%	&	4.283\%	&	4.757\%	&	5.086\%	&	5.248\%	\\	
	\hline
	\multicolumn{9}{c}{} \\ \hline
\st 	\multirow{15}{*}{\rotatebox{90}{70000}} & \multirow{5}{*}{\rotatebox{90}{0.075}}	&	5	&	0.189\%	&	0.214\%	&	0.399\%	&	0.481\%	&	0.521\%	&	0.559\%	\\	
\st	&	&	6	&	0.063\%	&	0.186\%	&	0.445\%	&	0.621\%	&	0.779\%	&	0.871\%	\\	
\st	&	&	7	&	0.047\%	&	0.268\%	&	0.497\%	&	0.623\%	&	0.749\%	&	0.840\%	\\	
\st	&	&	8	&	0.112\%	&	0.251\%	&	0.770\%	&	0.980\%	&	1.079\%	&	1.171\%	\\	
\st	&	&	9	&	0.314\%	&	0.584\%	&	1.449\%	&	1.771\%	&	2.038\%	&	2.179\%	\\	
	\cline{2-9}
\st &	\multirow{5}{*}{\rotatebox{90}{0.1}}
&	5	&	0.413\%	&	0.546\%	&	0.821\%	&	0.890\%	&	0.969\%	&	1.019\%	\\	
\st	&	&	6	&	0.257\%	&	0.504\%	&	0.992\%	&	1.280\%	&	1.476\%	&	1.565\%	\\	
\st	&	&	7	&	0.304\%	&	0.617\%	&	1.005\%	&	1.259\%	&	1.412\%	&	1.530\%	\\	
\st	&	&	8	&	0.281\%	&	0.701\%	&	1.563\%	&	1.794\%	&	1.980\%	&	2.113\%	\\	
\st	&	&	9	&	0.712\%	&	1.602\%	&	2.994\%	&	3.445\%	&	3.734\%	&	3.920\%	\\	
	\cline{2-9}
\st	&   \multirow{5}{*}{\rotatebox{90}{0.125}}
&	5	&	1.176\%	&	1.289\%	&	1.587\%	&	1.753\%	&	1.851\%	&	1.862\%	\\	
\st	&	&	6	&	0.607\%	&	1.246\%	&	2.117\%	&	2.580\%	&	2.830\%	&	3.003\%	\\	
\st	&	&	7	&	0.882\%	&	1.397\%	&	2.057\%	&	2.427\%	&	2.663\%	&	2.817\%	\\	
\st	&	&	8	&	0.863\%	&	1.999\%	&	3.221\%	&	3.487\%	&	3.742\%	&	3.875\%	\\	
\st	&	&	9	&	2.012\%	&	4.342\%	&	5.972\%	&	6.594\%	&	7.030\%	&	7.243\%	\\	
	\hline
 \multicolumn{9}{c}{} \\ \hline
  
\st 	\multirow{15}{*}{\rotatebox{90}{80000}} & \multirow{5}{*}{\rotatebox{90}{0.075}}	
	&	5	&	0.201\%	&	0.226\%	&	0.419\%	&	0.505\%	&	0.546\%	&	0.587\%	\\	
\st	&	&	6	&	0.068\%	&	0.198\%	&	0.473\%	&	0.660\%	&	0.827\%	&	0.925\%	\\	
\st	&	&	7	&	0.051\%	&	0.290\%	&	0.536\%	&	0.671\%	&	0.806\%	&	0.903\%	\\	
\st	&	&	8	&	0.125\%	&	0.276\%	&	0.843\%	&	1.071\%	&	1.178\%	&	1.278\%	\\	
\st	&	&	9	&	0.357\%	&	0.655\%	&	1.614\%	&	1.970\%	&	2.265\%	&	2.420\%	\\	\cline{2-9}

\st	& \multirow{5}{*}{\rotatebox{90}{0.1}}
&	5	&	0.450\%	&	0.588\%	&	0.877\%	&	0.949\%	&	1.034\%	&	1.086\%	\\	
\st	&	&	6	&	0.284\%	&	0.551\%	&	1.076\%	&	1.388\%	&	1.599\%	&	1.694\%	\\	
\st	&	&	7	&	0.344\%	&	0.687\%	&	1.110\%	&	1.387\%	&	1.554\%	&	1.683\%	\\	
\st	&	&	8	&	0.326\%	&	0.798\%	&	1.762\%	&	2.019\%	&	2.226\%	&	2.373\%	\\	
\st	&	&	9	&	0.850\%	&	1.874\%	&	3.460\%	&	3.970\%	&	4.298\%	&	4.510\%	\\	\cline{2-9}

\st	&   \multirow{5}{*}{\rotatebox{90}{0.125}}
&	5	&	1.347\%	&	1.441\%	&	1.752\%	&	1.931\%	&	2.037\%	&	2.047\%	\\	
\st	&	&	6	&	0.712\%	&	1.426\%	&	2.392\%	&	2.906\%	&	3.184\%	&	3.376\%	\\	
\st	&	&	7	&	1.065\%	&	1.638\%	&	2.377\%	&	2.796\%	&	3.064\%	&	3.238\%	\\	
\st	&	&	8	&	1.089\%	&	2.434\%	&	3.852\%	&	4.153\%	&	4.447\%	&	4.601\%	\\	
\st	&	&	9	&	2.684\%	&	5.535\%	&	7.440\%	&	8.173\%	&	8.691\%	&	8.944\%	\\	\hline
\end{tabular}
\caption{Network rebalancing without distribution costs. $\Delta Z^\circ_{Plex}$}\label{fig:Redis}
\end{table}

\begin{table}
\scriptsize
\centering
 \begin{tabular}{c|c|c|ccccc} \hline
  
  		\multicolumn{3}{c|}{Parameters}	& \multicolumn{5}{c}{$\gamma$}	\\ \hline	
\st $\phi$ &	$\alpha$	& K &	400.0	&	200.0	&	66.7	&	33.3	&	16.7		\\	\hline
\st 	\multirow{15}{*}{\rotatebox{90}{50000}} & \multirow{5}{*}{\rotatebox{90}{0.075}}	
&	5	&	0.130\%	&	0.160\%	&	0.426\%	&	1.613\%	&	0.203\%	\\
\st & &	6	&	0.028\%	&	0.152\%	&	0.363\%	&	0.891\%	&	1.015\%		\\
\st & &	7	&	0.012\%	&	0.192\%	&	0.461\%	&	1.011\%	&	2.318\%	\\
\st & &	8	&	0.100\%	&	0.266\%	&	0.875\%	&	1.692\%	&	2.433\%		\\
\st & &	9	&	0.248\%	&	0.552\%	&	1.378\%	&	1.788\%	&	2.225\%	\\
\cline{2-8}

\st &	\multirow{5}{*}{\rotatebox{90}{0.1}}	
&	5	&	0.294\%	&	0.464\%	&	1.046\%	&	3.078\%	&	0.436\%	\\
\st & &	6	&	0.190\%	&	0.407\%	&	0.855\%	&	1.737\%	&	1.905\%		\\
\st & &	7	&	0.203\%	&	0.464\%	&	0.951\%	&	1.943\%	&	4.233\%		\\
\st & &	8	&	0.225\%	&	0.763\%	&	1.722\%	&	3.201\%	&	4.291\%	\\
\st & &	9	&	0.534\%	&	1.449\%	&	2.625\%	&	3.449\%	&	4.106\%		\\
\cline{2-8}

	\st &	\multirow{5}{*}{\rotatebox{90}{0.125}}	
&	5	&	0.913\%	&	1.212\%	&	2.008\%	&	6.167\%	&	0.851\%	\\
\st & &	6	&	0.482\%	&	0.966\%	&	1.718\%	&	3.331\%	&	3.471\%	\\
\st & &	7	&	0.561\%	&	1.052\%	&	1.860\%	&	3.810\%	&	7.612\%	\\
\st & &	8	&	0.559\%	&	1.782\%	&	3.267\%	&	5.979\%	&	7.493\%	\\
\st & &	9	&	1.293\%	&	3.171\%	&	4.957\%	&	6.329\%	&	7.126\%		\\
	\hline

	\multicolumn{8}{c}{} \\ \hline
\st 	\multirow{15}{*}{\rotatebox{90}{70000}} & \multirow{5}{*}{\rotatebox{90}{0.075}}	
&	5	&	0.146\%	&	0.178\%	&	0.468\%	&	1.772\%	&	0.223\%	\\
\st & &	6	&	0.032\%	&	0.171\%	&	0.407\%	&	0.999\%	&	1.137\%		\\
\st & &	7	&	0.014\%	&	0.222\%	&	0.529\%	&	1.159\%	&	2.659\%	\\
\st & &	8	&	0.121\%	&	0.316\%	&	1.030\%	&	1.987\%	&	2.859\%		\\
\st & &	9	&	0.310\%	&	0.675\%	&	1.666\%	&	2.158\%	&	2.682\%	\\
\cline{2-8}

\st &	\multirow{5}{*}{\rotatebox{90}{0.1}}
&	5	&	0.344\%	&	0.531\%	&	1.183\%	&	3.483\%	&	0.491\%		\\
\st & &	6	&	0.228\%	&	0.477\%	&	0.993\%	&	2.013\%	&	2.206\%		\\
\st & &	7	&	0.251\%	&	0.560\%	&	1.133\%	&	2.314\%	&	5.062\%	\\
\st & &	8	&	0.288\%	&	0.954\%	&	2.123\%	&	3.938\%	&	5.286\%		\\
\st & &	9	&	0.713\%	&	1.881\%	&	3.352\%	&	4.394\%	&	5.223\%		\\
\cline{2-8}

\st	&   \multirow{5}{*}{\rotatebox{90}{0.125}}
&	5	&	1.150\%	&	1.473\%	&	2.398\%	&	7.393\%	&	1.008\%	\\
\st & &	6	&	0.628\%	&	1.215\%	&	2.120\%	&	4.109\%	&	4.279\%		\\
\st & &	7	&	0.759\%	&	1.370\%	&	2.375\%	&	4.867\%	&	9.835\%		\\
\st & &	8	&	0.799\%	&	2.443\%	&	4.378\%	&	7.998\%	&	10.073\%	\\
\st & &	9	&	1.963\%	&	4.595\%	&	6.983\%	&	8.899\%	&	9.998\%	\\
\hline

\multicolumn{8}{c}{} \\ \hline
  
\st 	\multirow{15}{*}{\rotatebox{90}{80000}} & \multirow{5}{*}{\rotatebox{90}{0.075}}	
&	5	&	0.156\%	&	0.188\%	&	0.492\%	&	1.864\%	&	0.234\%		\\
\st & &	6	&	0.035\%	&	0.184\%	&	0.434\%	&	1.063\%	&	1.210\%		\\
\st & &	7	&	0.016\%	&	0.241\%	&	0.571\%	&	1.250\%	&	2.870\%		\\
\st & &	8	&	0.135\%	&	0.349\%	&	1.130\%	&	2.177\%	&	3.132\%		\\
\st & &	9	&	0.353\%	&	0.759\%	&	1.861\%	&	2.407\%	&	2.988\%		\\
\cline{2-8}

\st	& \multirow{5}{*}{\rotatebox{90}{0.1}}
&	5	&	0.375\%	&	0.572\%	&	1.266\%	&	3.729\%	&	0.524\%		\\
\st & &	6	&	0.253\%	&	0.523\%	&	1.079\%	&	2.187\%	&	2.395\%		\\
\st & &	7	&	0.284\%	&	0.625\%	&	1.253\%	&	2.558\%	&	5.610\%		\\
\st & &	8	&	0.335\%	&	1.090\%	&	2.402\%	&	4.450\%	&	5.979\%	\\
\st & &	9	&	0.856\%	&	2.210\%	&	3.890\%	&	5.091\%	&	6.045\%		\\
\cline{2-8}

\st	&   \multirow{5}{*}{\rotatebox{90}{0.125}}
&	5	&	1.322\%	&	1.650\%	&	2.656\%	&	8.208\%	&	1.110\%		\\
\st & &	6	&	0.741\%	&	1.395\%	&	2.402\%	&	4.653\%	&	4.841\%	\\
\st & &	7	&	0.920\%	&	1.614\%	&	2.757\%	&	5.651\%	&	11.516\%	\\
\st & &	8	&	1.017\%	&	3.000\%	&	5.274\%	&	9.622\%	&	12.169\%	\\
\st & &	9	&	2.649\%	&	5.926\%	&	8.777\%	&	11.166\%	&	12.522\%	\\
\hline
\end{tabular}
\caption{Network rebalancing with distribution costs, $\varrho=8.30$ cts/km $\times$ Ton}\label{fig:RedisDist}
\end{table}

\begin{table}
	\scriptsize
	\centering
	\begin{tabular}{c|c|cccccc} \hline
		\multicolumn{2}{c|}{Parameters}	& \multicolumn{5}{c}{$\gamma$}	\\	\hline
		\st 	$\alpha$	& K &	400.0	&	200.0	&	66.7	&	33.3	&	16.7	&	8.3	\\	\hline
		\st 	\multirow{5}{*}{\rotatebox{90}{0.075}}	
	&	5	&	-0.45	\%&	-0.79	\%&	-0.93	\%&	-0.90	\%&	-0.90	\%&	-0.87	\%	\\
    \st		&	6	&	-0.51	\%&	-0.70	\%&	-1.27	\%&	-1.30	\%&	-1.48	\%&	-1.43	\%	\\
    \st		&	7	&	-0.69	\%&	-0.77	\%&	-1.09	\%&	-1.25	\%&	-1.40	\%&	-1.45	\%	\\
    \st		&	8	&	-0.18	\%&	-0.87	\%&	-1.82	\%&	-1.90	\%&	-1.90	\%&	-1.91	\%	\\
    \st		&	9	&	-0.66	\%&	-2.03	\%&	-3.22	\%&	-3.45	\%&	-3.43	\%&	-3.39	\%	\\
	 \cline{1-8}
		
		\st 	\multirow{5}{*}{\rotatebox{90}{0.1}}	
			&	5	&	-0.85	\%&	-0.95	\%&	-1.00	\%&	-0.97	\%&	-0.90	\%&	-0.88	\%	\\
		\st  	&	6	&	-0.51	\%&	-0.86	\%&	-1.45	\%&	-1.59	\%&	-1.48	\%&	-1.46	\%	\\
		\st 	&	7	&	-0.86	\%&	-0.99	\%&	-1.36	\%&	-1.39	\%&	-1.49	\%&	-1.48	\%	\\
		\st  	&	8	&	-0.63	\%&	-1.14	\%&	-1.97	\%&	-1.89	\%&	-1.90	\%&	-1.91	\%	\\
		\st  	&	9	&	-1.10	\%&	-2.58	\%&	-3.49	\%&	-3.59	\%&	-3.50	\%&	-3.39	\%	\\			
		\cline{1-8}
		
		\st 	\multirow{5}{*}{\rotatebox{90}{0.125}}	
			&	5	&	-1.10	\%&	-1.32	\%&	-1.11	\%&	-0.97	\%&	-0.91	\%&	-0.89	\%	\\
		\st  	&	6	&	-0.89	\%&	-1.37	\%&	-1.55	\%&	-1.60	\%&	-1.53	\%&	-1.51	\%	\\
		\st  	&	7	&	-0.86	\%&	-1.27	\%&	-1.44	\%&	-1.62	\%&	-1.52	\%&	-1.48	\%	\\
		\st  	&	8	&	-1.01	\%&	-2.04	\%&	-2.08	\%&	-2.01	\%&	-1.98	\%&	-1.91	\%	\\
		\st  	&	9	&	-2.71	\%&	-3.45	\%&	-3.83	\%&	-3.66	\%&	-3.50	\%&	-3.39	\%	\\			
		\hline
	\end{tabular}
	\caption{ Change in complexity costs ($\Delta C_\alpha$). Network rebalancing without distribution costs.}\label{fig:RedisComp}
\end{table}

\begin{table}
	\scriptsize
	\centering
	\begin{tabular}{c|c|ccccc} \hline
		
		 \multicolumn{2}{c|}{Parameters}	& \multicolumn{5}{c}{$\gamma$}	\\ \hline	
		\st 	$\alpha$	& K &	400.0	&	200.0	&	66.7	&	33.3	&	16.7		\\	\hline
		\st 	 \multirow{5}{*}{\rotatebox{90}{0.075}}	
				&5	&	-0.50	\%&	-0.74	\%&	-1.37	\%&	-3.41	\%&	-0.75	\%	\\
		\st 	&	6	&	-0.54	\%&	-0.49	\%&	-1.34	\%&	-2.16	\%&	-2.27	\%	\\
		\st  	&	7	&	-0.53	\%&	-0.72	\%&	-1.41	\%&	-2.41	\%&	-4.48	\%	\\
		\st  	&	8	&	-0.17	\%&	-0.92	\%&	-2.25	\%&	-4.10	\%&	-4.62	\%	\\
		\st  	&	9	&	-0.65	\%&	-2.09	\%&	-3.31	\%&	-4.61	\%&	-4.46	\%	\\
		\cline{1-7}
		
		\st 	\multirow{5}{*}{\rotatebox{90}{0.1}}	
				&5	&	-0.94	\%&	-1.10	\%&	-1.41	\%&	-3.50	\%&	-0.75	\%	\\
		\st  	&	6	&	-0.54	\%&	-0.88	\%&	-1.60	\%&	-2.23	\%&	-2.27	\%	\\
		\st  	&	7	&	-0.76	\%&	-1.16	\%&	-1.54	\%&	-2.62	\%&	-4.61	\%	\\
		\st  	&	8	&	-0.49	\%&	-1.72	\%&	-2.60	\%&	-4.30	\%&	-4.62	\%	\\
		\st  	&	9	&	-0.98	\%&	-3.17	\%&	-3.96	\%&	-4.86	\%&	-4.76	\%	\\
		\cline{1-7}
		
		\st 	\multirow{5}{*}{\rotatebox{90}{0.125}}	
			&5	&	-1.21	\%&	-1.49	\%&	-1.51	\%&	-3.46	\%&	-0.76	\%	\\
		\st  	&	6	&	-0.54	\%&	-1.45	\%&	-1.73	\%&	-2.40	\%&	-2.27	\%	\\
		\st  	&	7	&	-0.76	\%&	-1.31	\%&	-1.85	\%&	-2.94	\%&	-4.64	\%	\\
		\st  	&	8	&	-0.89	\%&	-2.40	\%&	-2.78	\%&	-4.44	\%&	-4.62	\%	\\
		\st  	&	9	&	-2.59	\%&	-3.85	\%&	-4.10	\%&	-5.02	\%&	-4.76	\%	\\
		\hline
	\end{tabular}
	\caption{Change in complexity costs ($\Delta C_\alpha$). Network rebalancing, $\varrho=8.30$ cts/km $\times$ Ton.}\label{fig:RedisDistComp}
\end{table}

\subsection{Network rationalisation}

Experiments conducted to assess the impact of the \textit{network rationalisation} strategy returned positive improvements only in a very limited number of instances. Improvements in profitability were obtained only for high values of parameter $\alpha$, for a small number of facilities, and for large transport costs. 

For example, for the case without distribution costs, the best improvements were attained for cases with three facilities, $\phi=80k$, and transport costs off $400cts/km \times ton$. For $\alpha=0.1$ the maximum improvement was $0.20\%$, whereas for an $\alpha=0.125$ this value increased to $3.55\%$. In overall, in 86.8\% of the cases it was not possible to attain any improvement in the network's profitability.

Taking into account distribution costs, and considering a value of $0.1$ for parameter $\alpha$, the largest improvements observed were always obtained for of $\varrho = 200$ and transport costs equal to $400$. These values were $0.275\%$ for $\phi=50k$; $0.325\%$ for $\phi=70k$; and $0.358\%$ for $\phi=80k$. In all those cases, the best improvement was obtained for networks with $K=3$ facilities. For $\alpha=0.125$, the procedure returns the following maximum improvements for parameters $\varrho=200$ and $\gamma=400$: $5.31\%$ for $\phi=50k$; $7.35\%$ for $\phi=70k$; and $9.11\%$ for $\phi=80k$. In all cases, the network size for which the largest improvement was attained was $K=3$.  Notice that improvement was  only attained for small networks, where no rationalisation is indeed required. For over-dimensioned networks, no improvement was found.  Indeed, no improvement in profitability was obtained in 82.3\% of the instances. 

Regarding the cost of complexity, the elimination of nodes allows the system to attain some savings in $C_\alpha$ for small networks. These savings fall in a range between 1\% and 25\%. There is no evidence that any reduction in the cost of complexity could be obtained with this strategy.

\subsection{Network reduction}

In this section we present the results of the numerical experiments conducted for assessing the \textit{network reduction} strategy. We evaluate two alternative approaches, the first one removes the demand nodes associated with the eliminated facility; while in the second one the demand nodes are reassigned to facilities that remain open. 

This strategy renders better results than the other two approaches. While the maximum increase in profit attained with the \textit{network rebalancing} strategy --without considering distribution costs- was 8.94\%, the highest improvement with the \textit{network reduction} strategy is 24.2\%. Likewise, when distribution costs are considered, \textit{network reduction} can bring improvements up to 120\% (considering only values of $\alpha$ below 0.125), whereas the maximal improvement brought by the \textit{network rebalancing} strategy (for the same values of $\alpha$) is just above 6\%.

Without distribution costs, this strategy attains reductions in complexity costs  to 12.7\% when the demand nodes left uncovered are abandoned. If, instead, uncovered demand nodes are reallocated to open facilities, the total cost of complexity grows between 1\% and 6.1\%. When distribution costs between the central facility and the regional facilities are introduced, the \textit{network reduction} strategy can bring savings around 26\% if uncovered nodes are abandoned; otherwise, if demand is reallocated, the cost of complexity increases up to 11\%. These results confirm our conjecture that significant savings in complexity costs can only be attained with strategies that involve abandoning unprofitable markets.

It is important to highlight that in our model operation costs are assumed fixed for each facility. This implies that the natural increase in operating expenses, resulting from the necessary increase in a facility's capacity to serve a larger market, is ignored. Consequently, the change in profit reported in the tables should be seen as an upper bound on the potential increase that would result from reallocating demand.

Table \ref{tab:Reduc} shows the results obtained for a subset of our experiments for the case without distribution costs. Results are reported for four different values of parameter $\gamma$ (transport cost) and three different values of $\alpha$. The table reports the initial and final number of facilities, $K_0$ and $K$. It also includes the value of the objective function, $Z^\circ_{Plex}(\mathcal N^S, S)$, and the associated cost of complexity, $C_\alpha (\mathcal N^S, \mc S)$, evaluated in the solution of the $K$-Median problem. It also includes the relative improvement, $\Delta Z^\circ_{Plex}(\mathcal N'^{S'}, \mc S')$, in the objective function for the reduced network, $\mathcal N'^{S'}$ , and the corresponding change in the cost of complexity, $\Delta C_\alpha (\mathcal N'^{S'}, \mc S')$. Finally, values for $\Delta Z^\circ_{Plex}(\mathcal N^{S'} \mc S')$ and $\Delta C_\alpha (\mathcal N^{S'}, \mc S')$, are presented for the network with reallocated demand nodes, $\mathcal N^{S'}$.

All the reported values were obtained for network configurations with parameter $\phi=80k$. The symbol $\leq$ indicates that no improvement was obtained for instances where the initial number of facilities was equal or below the given value. 

A few regularities can be observed in the table. In general, larger improvements in profit are obtained when the uncovered demand nodes are reallocated to other facilities. However, savings in complexity costs are only attained when those demand nodes are abandoned. Indeed, after reallocating abandoned nodes, the complexity costs increase substantially. The positive impact of \textit{network reduction} is lager for higher transport costs, this is also the case for the savings in complexity costs.

Table \ref{tab:ReducDist} presents the same information as Table \ref{tab:Reduc} when distribution costs are considered. For these cases, the network's revenue is calculated using equation \eqref{eq:revalt}. Please notice that results in the north-east corner of the table, corresponding to values of $\alpha=0.125$ and high distribution and transport costs ($\varrho= 200$, $\gamma =400$), must be taken with a caveat. In those cases, the network seems to be extremely over dimensioned, and the combination of high distribution and high complexity cost renders the original network highly inefficient. This implies that large reductions in the network size may have a huge impact in profits. However, situations like this will hardly appear practice, and therefore the results are only presented for completeness.

In general, the impact of \textit{network reduction} with distribution costs appears to be higher than without them. Other regularities observed in the previous cases, are also present in Table \ref{tab:ReducDist}: higher benefits are obtained when the demand nodes are reassigned among facilities that remain open but, again, there is a significant increase in the costs of complexity. Finally, the higher the distribution costs, the larger the impact on profit and complexity-derived costs of the \textit{network reduction} strategy. 

Before concluding this section, it is important mentioning that, although we have limited the results presented here to cases where the value of parameter $\phi$ is set to $80k$, a larger set of experiments was ran for values of $50k$ and $70k$. In those cases, the results presented the same regularities observed in the tables in this section, but with significantly smaller improvements. We therefore decided to leave them out of the report.

\begin{table}[!h]
\footnotesize\centering
 \begin{tabular}{cccccccc} \hline
  \st	$\gamma$	&	400	&	400	&	400	& 	400	&	400	&	400 & 400	\\	\hline
\st	$\alpha$	&	0.075	&	0.1	&	0.1	&	0.1	&	0.125	&	0.125 & 0.125	\\	\hline
\st	$K_0=\abs{S}$	&$\leq$	9	&	9	&	8	&	$\leq$7&	9	&	8 &	$\leq$7 \\	
\st	$K=\abs{S'}$	&	N.I.&	7	&	7	&	N.I.	&	6	&	6	& N.I.\\	
\st	$Z^\circ_{Plex}(\mathcal N^S, S)$	&		&	462900.4	&	503282.0	&		&	269293.4	&	305117.2	& \\	
\st $C_\alpha (\mathcal N^S, \mc S)$  	& &774428.1 &		792659.5 &  &		968035.2&		990824.3& \\
\st $\Delta Z^\circ_{Plex} (\mathcal N'^{S'}, \mc S')$ 	&		&	1.3 \% 	&	0.2 \%	&		&	7.7 \% &	2.9 \%	&\\	
\st	$\Delta C_\alpha (\mathcal N'^{S'}, \mc S')$ 	&		&	-7.3 \% 	&	-3.9 \%	&		&	-12.7 \% &	-9.2 \%	&\\	
\st $\Delta Z^\circ_{Plex} (\mathcal N^{S'}, \mc S')$		&		&	15.5 \%	&	6.3 \%	&		&	24.2 \%	&	9.6 \%	&\\	
\st $\Delta C_\alpha (\mathcal N^{S'}, \mc S')$		&		&	3.4 \%	&	1.0 \%	&		&	4.8 \%	&	2.4 \%	&\\	
		&		&		&		&		&		&	&	\\ 
		
\st	$\gamma$	&	200	&	200	&	200	&	200	&	200	&	200 & 200	\\	\hline
\st	$\alpha$	&	0.075	& 0.1  &	0.1	&	0.1	&	0.125	&	0.125	&	0.125	\\	\hline
\st	$K_0=\abs{S}$	&	$\leq$9	&	9	&	$\leq$8	&	& 9	&	8	&	$\leq$7	\\	
\st	$K=\abs{S'}$	&	N.I. 	&	8	&	N.I. &	&	7	&	7	&	N.I.	\\	
\st	$Z^\circ_{Plex}(\mathcal N^{S},S)$	&		&	532048.2	&	& 	&	327497.6	&	367294.6	&		\\	
\st $C_\alpha (\mathcal N^S, \mc S)$  	& & 818202.1 & & &		1022752.6 &	1052262.1 &  \\
\st	$\Delta Z^\circ_{Plex} (\mathcal N'^{S'}, \mc S')$	&		&	0.7 \% 	&		&	& 4.5 \% 	&	1.2 \%	 &		\\	
\st	$\Delta C_\alpha (\mathcal N'^{S'}, \mc S')$ 	&		& -3.3\% &  & & -7.2\% &-3.9\% & \\ 
\st	$\Delta Z^\circ_{Plex} (\mathcal N^{S'}, \mc S')$	&		&	8.6 \%	&		&	& 23.6 \%	&	10.2 \%	&		\\	
\st$\Delta C_\alpha (\mathcal N^{S'}, \mc S')$		& & 2.9\% & & & 5.1\% & 2.1\%  & \\
	&	&		&		&		&		&		&		\\	
	
\st	$\gamma$	&	100	&	100 & 100	&	100	&	100	&	100	&	100	\\	\hline
\st	$\alpha$	&	0.075	&	0.1	&	0.1 & 0.1	&	0.125	&	0.125	&	0.125	\\	\hline
\st	$K_0=\abs{S}$	&	$\leq$9	&	9	&	$\leq$8	&	& 9	&	8	&	$\leq$7	\\	
\st	$K=\abs{S'}$	&	N.I.	&	8	&	N.I.	& &	7	&	7	&	N.I.	\\	
\st	$Z^\circ_{Plex}(\mathcal N^{S},S)$	&		&	566622.0	&		&	& 356599.8	&	398383.3	&		\\	
\st $C_\alpha (\mathcal N^S, \mc S)$  	& & 840089.1 &  & &1050111.4 &	1082981.0  &  \\
\st	$\Delta Z^\circ_{Plex} (\mathcal N'^{S'}, \mc S')$	&		&	0.5 \%	 &		&	& 3.3 \% 	&	0.6 \%	&		\\
\st	$\Delta C_\alpha (\mathcal N'^{S'}, \mc S')$ 	& &-3.2\% & & &-7.2\% & -3.9\% & \\
\st	$\Delta Z^\circ_{Plex} (\mathcal N^{S'}, \mc S')$	&		&	8.5 \%	&		&	& 23.3 \%	&	10.4 \%	&		\\	
\st $\Delta C_\alpha (\mathcal N^{S'}, \mc S')$		& & 3.1\% & & & 5.9\% & 2.6\%  & \\
	&	&		&		&		&		&		&		\\	

\st	$\gamma$	&	66.6	&	66.6 & 66.6	&	66.6	&	66.6	&	66.6	&	66.6	\\	\hline
\st	$\alpha$	&	0.075	&	0.1	&	0.1	&	0.1 & 0.125	&	0.125	&	0.125	\\	\hline
\st	$K_0=\abs{S}$	&	$\leq$ 9	&	9	&	$\leq$ 8 &	&	9	&	8	&	$\leq$7	\\	
\st	$K=\abs{S'}$	&	N.I.	&	8	&	N.I. & 	&	7	&	7	&	N.I.	\\	
\st	$Z^\circ_{Plex}(\mathcal N^{S},S)$	&		&	578146.7	&		&	& 366300.5	&	408746.3	&		\\	
\st $C_\alpha (\mathcal N^S, \mc S)$  	& & 847384.8& & & 1059230.9 & 1093220.7 & \\
\st	$\Delta Z^\circ_{Plex} (\mathcal N'^{S'}, \mc S')$	&		&	0.4 \%	&		&	& 3.0 \%	&	0.4 \%	&		\\	
\st	$\Delta C_\alpha (\mathcal N'^{S'}, \mc S')$ 	& &-3.2\% & & &-7.2\% & -3.9\% & \\
\st	$\Delta Z^\circ_{Plex} (\mathcal N^{S'}, \mc S')$	&		&	8.5 \%	&		&	& 23.2 \%	&	10.4 \%	&		\\
\st $\Delta C_\alpha (\mathcal N^{S'}, \mc S')$	& &3.2\% & & & 6.1\%&2.8\% & \\	
	&	&		&		&		&		&		&		\\	\hline
 \end{tabular}
\caption{Network reduction with $\phi=80k$.  $C_p(\mathcal N^{\mc S}, \mc S)=6.96$. (N.I.: No improvement.)}\label{tab:Reduc}
\end{table}

\newpage

\begin{landscape}
\begin{table}
\vspace{-1cm}
\scriptsize
\begin{tabular}{c|cccccccccccc}\\
\st $\gamma$	&	400	&	400	&	400	&	400	&	400	&	400	&	400	&	400	&\itshape	400	&\itshape	400	&\itshape	400	&	\itshape 400	\\	\hline
\st $\varrho$	&	200	&	200	&	200	&	200	&	200	&	200	&	200	&	200	&	\itshape 200	&	\itshape 200	&	\itshape 200	&	\itshape 200	\\	\hline
\st $\alpha$	&	0.075	&	0.075	&	0.075	&	0.075	&	0.1	&	0.1	&	0.1	&	0.1	&	\itshape0.125	&\itshape	0.125	&	\itshape0.125	&	\itshape0.125	\\	\hline
\st $K_0=\abs{S}$	&	9	&	8	&	7	&	6	&	9	&	8	&	7	&	6	&	\itshape9	&\itshape	8	&	\itshape7	&	\itshape 6	\\	
\st $K=\abs{S'}$	&	5	&	5	&	5	&	5	&	4	&	4	&	4	&	5	&\itshape	2	&	\itshape3	&\itshape3	&\itshape	5	\\	
\st $Z^\circ_{Plex}(\mathcal N^{S},S)$	&	269783.3	&	315652.5	&	356805.9	&	409682.8	&	124102.3	&	164459.3	&	200728.1	&	254898.7	&	\itshape-21578.7	&	\itshape13266.0	&	\itshape44650.3	&\itshape	100114.5	\\	
\st $C_\alpha (\mathcal N^S, \mc S)$	&	437043.1 &	453579.7	&468233.5	&464352.4	&	582724.1	&604772.9	&624311.3	&619136.6		&\itshape728405.1	&\itshape755966.1& \itshape	780389.1	&\itshape773920.7 	\\	
\st $\Delta Z^\circ_{Plex} (\mathcal N'^{S'}, \mc S')$	&	16.6 \%	&	14.1 \%	 &	10.8 \%	&	7.3 \%	&	63.8 \%	&	43.0 \% &	28.7 \%	&	13.4 \%	 &	\itshape621.1 \% &	\itshape824.9 \% 	&\itshape	180.1 \%	&	\itshape38.6 \% 	\\	
\st $\Delta C_\alpha (\mathcal N'^{S'}, \mc S')$	&-19.1\%&	-12.9\%&	-7.9\%&	-2.9\%&	-26.1\%&	-19.7\%&	-14.5\%&	-2.9\%&		\itshape-41.5\%&	\itshape-27.3\%&	\itshape-24.2\%&	\itshape-2.9\% 	\\	
\st	$\Delta Z^\circ_{Plex} (\mathcal N^{S'}, \mc S')$& 	64.5 \%	&	37.7 \%	&	21.8 \%	&	9.1 \%	&	118.2 \%	&	57.9 \%	&	29.3 \%	&	13.5 \%	&	\itshape510.2 \%	&\itshape	612.8 \%	&\itshape	111.8 \%	&\itshape	31.5 \%	\\	
\st $\Delta C_\alpha (\mathcal N^{S'}, \mc S')$ & 8.2\%&	5.9\%&	2.6\%&	1.8\%&	9.7\%&	7.2\%&	3.9\%&	1.8\%&	\itshape22.8\%&	\itshape11.4\%&	\itshape7.9\%&	\itshape1.8\% \\

	&		&		&		&		&		&		&		&		&		&		&		&		\\	
\st $\gamma$	&	200	&	200	&	200	&	200	&	200	&	200	&	200	&	200	&	200	&	200	&	200	&	200	\\	\hline
\st $\varrho$	&	100	&	100	&	100	&	100	&	100	&	100	&	100	&	100	&	100	&	100	&	100	&	100	\\	\hline
\st $\alpha$	&	0.075	&	0.075	&	0.075	&	0.075	&	0.1	&	0.1	&	0.1	&	0.1	&	0.125	&	0.125	&	0.125	&	0.125	\\	\hline
\st $K_0=\abs{S}$	&	9	&	8	&	7	&	6	&	9	&	8	&	7	&	6	&	9	&	8	&	7	&	6	\\	
\st $K=\abs{S'}$	&	8	&	7	&	6	&	5	&	8	&	7	&	6	&	5	&	6	&	5	&	5	&	5	\\	
\st $Z^\circ_{Plex}(\mathcal N^{S},S)$	&	539424.7	&	586677.3	&	630653.0	&	686976.4	&	357566.6	&	396835.7	&	432835.4	&	486434.0	&	175708.4	&	206994.1	&	235017.8	&	285891.5	\\	
\st $C_\alpha (\mathcal N^S, \mc S)$	&	545574.5&	569524.8&	593452.7&	601627.3&	727432.6&	759366.4&	791270.3&	802169.7&	909290.8&	949207.9&	989087.9&	1002712.1 	\\	
\st $\Delta Z^\circ_{Plex} (\mathcal N'^{S'}, \mc S')$	&	1.2 \%	 &	1.1 \% 	&	1.0 \%	 &	0.9 \% &	3.6 \%	 &	3.2 \%&	3.0 \%	 &	2.6 \%	&	15.4 \%	&	13.1 \%&	8.6 \%	&	6.8 \%	 \\	
\st $\Delta C_\alpha (\mathcal N'^{S'}, \mc S')$	&	-3.6\%&	-3.4\%&	-3.3\%&	-3.3\%&	-3.6\%&	-3.4\%&	-3.3\%&	-3.3\%&	-14.6\%& -14.0\%& -9.0\%& -3.3\%  \\	
\st	$\Delta Z^\circ_{Plex} (\mathcal N^{S'}, \mc S')$ &	9.9 \%	&	6.9 \%	&	6.3 \%	&	6.0 \%	&	13.6 \%	&	8.3 \%	&	7.3 \%	&	6.9 \%	&	61.8 \%	&	42.2 \%	&	25.2 \%	&	9.2 \%	\\	
\st $\Delta C_\alpha (\mathcal N^{S'}, \mc S')$&2.6\%&	4.0\%&	3.9\%&	3.7\%&	2.6\%&	4.0\%&	3.9\%&	3.7\%&	9.8\%&	11.1\%&	6.6\%&	3.7\% \\
	&		&		&		&		&		&		&		&		&		&		&		&		\\	
\st $\gamma$	&	200	&	& & 	&	200	&	200	& & &	200	&	200	&			&		\\	\hline
\st $\varrho$	&	66.6	&	&	& & 66.6	&	66.6 &  &	&	66.6	&	66.6	&		\\	\hline
\st $\alpha$	&	0.075	&	& &	&	0.1	&	0.1	& & &	0.125	&	0.125	&		&			\\	\hline
\st $K_0=\abs{S}$	&	9	&	& &	&	9	&	8	& & &	9	&	8	&		&	\\	
\st $K=\abs{S'}$	&	8	&	& &	&	7	&	7 & &	&	7	&	7	&		&	\\	
\st $Z^\circ_{Plex}(\mathcal N^{S},S)$	&	600359.9	&	& &	&	411235.5	&	450944.6	& & &	222111.1	&	254068.8	&		&			\\	
\st $C_\alpha (\mathcal N^S, \mc S)$	&	567373.2	&	&  & 	&	756497.6	&	787503.3	& & &	945622.0	&	984379.1	&		&		\\	
\st $\Delta Z^\circ_{Plex} (\mathcal N'^{S'}, \mc S')$	&	0.2 \%	&	&  & 	&	2.4 \%	 &	0.4 \%	& & &	11.0 \%	&	3.7 \%	&		&	\\	
\st $\Delta C_\alpha (\mathcal N'^{S'}, \mc S')$	&	-3.7\%	& & &  		&	-7.7\%	&	-3.9\% & &	&	-7.7\%	&	-3.9\%	&		&	\\	
\st 	$\Delta Z^\circ_{Plex} (\mathcal N^{S'}, \mc S')$	&	8.3 \%	&	& &	&	21.7 \%	&	9.8 \%	& & &	35.6 \%	&	15.2 \%	&		&	\\	
\st $\Delta C_\alpha (\mathcal N^{S'}, \mc S')$&	2.4\%	&	& &	&	5.3\%	&	2.7\%	& & &	5.3\%	&	2.7\%	&		&	\\
\hline
 \end{tabular}
\caption{Network reduction with $\phi=80k$.  $C_p(\mathcal N^{\mc S}, \mc S)=6.96$.}\label{tab:ReducDist}
\end{table}
\end{landscape}

\subsection{Illustrative example}

In this section we present the results of applying the \textit{network rebalancing} and \textit{network reduction} strategies in an instance with $\alpha=0.125$, $\gamma= 200 cts/km \times ton$, and $\varrho=66.6 cts/km \times ton$. 

\begin{figure}
	\begin{center}
		\begin{tabular}{cc}
			\includegraphics[width=6cm]{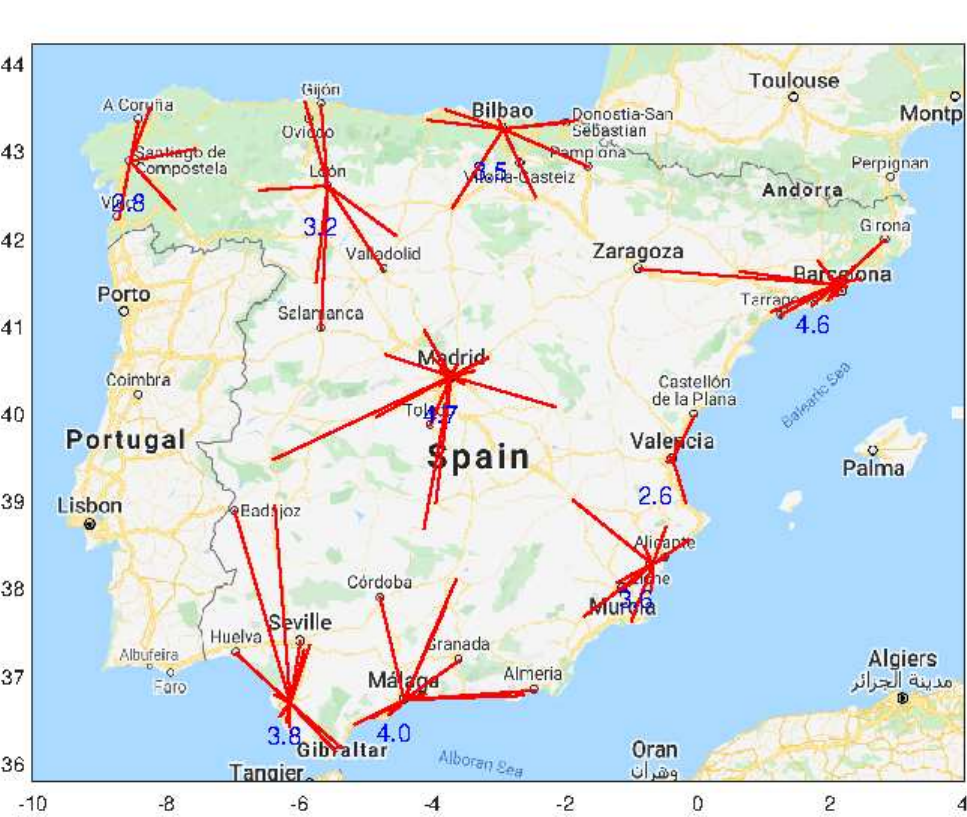} &
			\includegraphics[width=6cm]{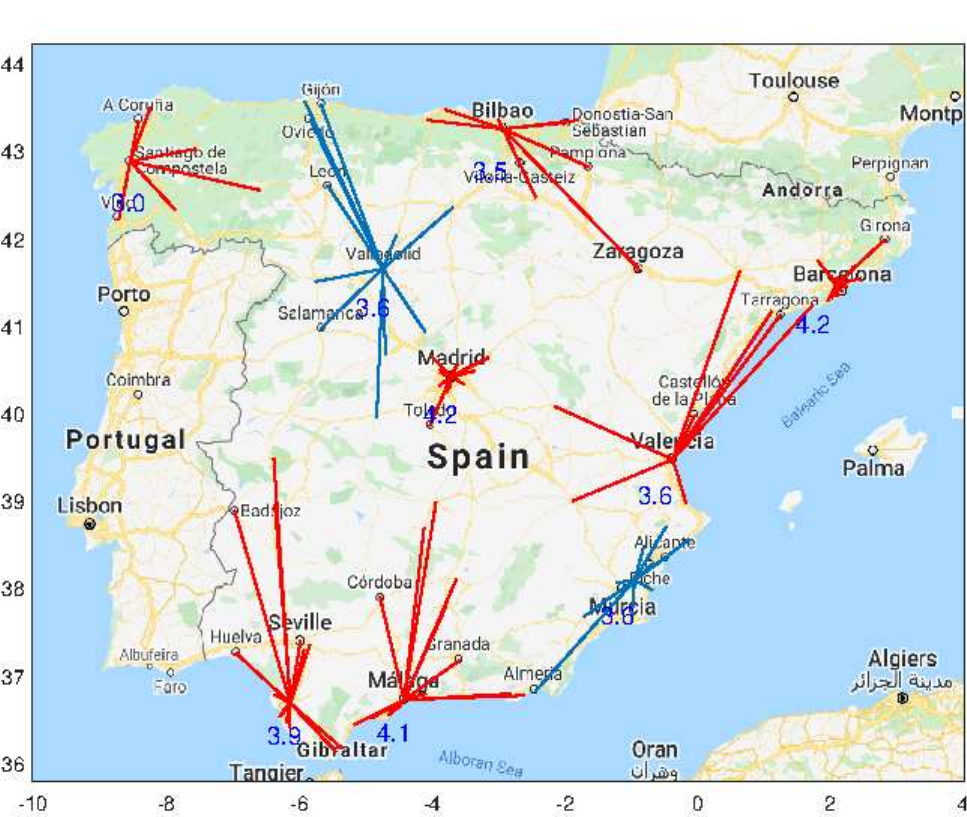} \\ (a) $Z^K$/$Z^K_{Plex}$ & (b) Rebalancing\\
			\includegraphics[width=6cm]{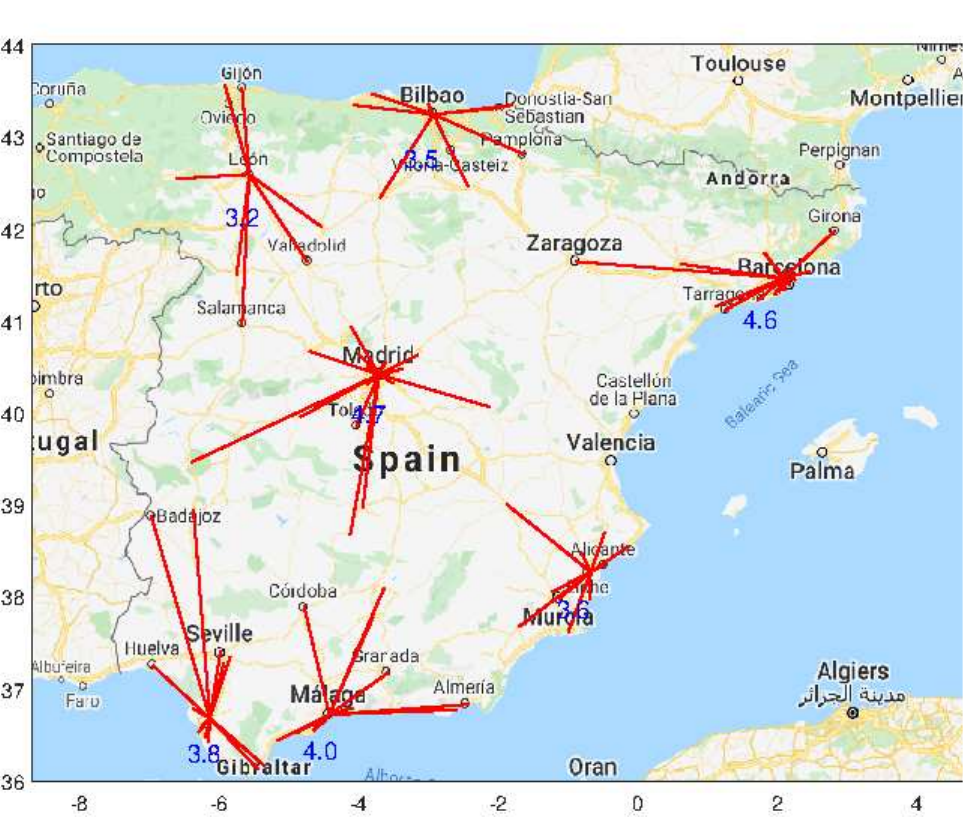}&
			\includegraphics[width=6cm]{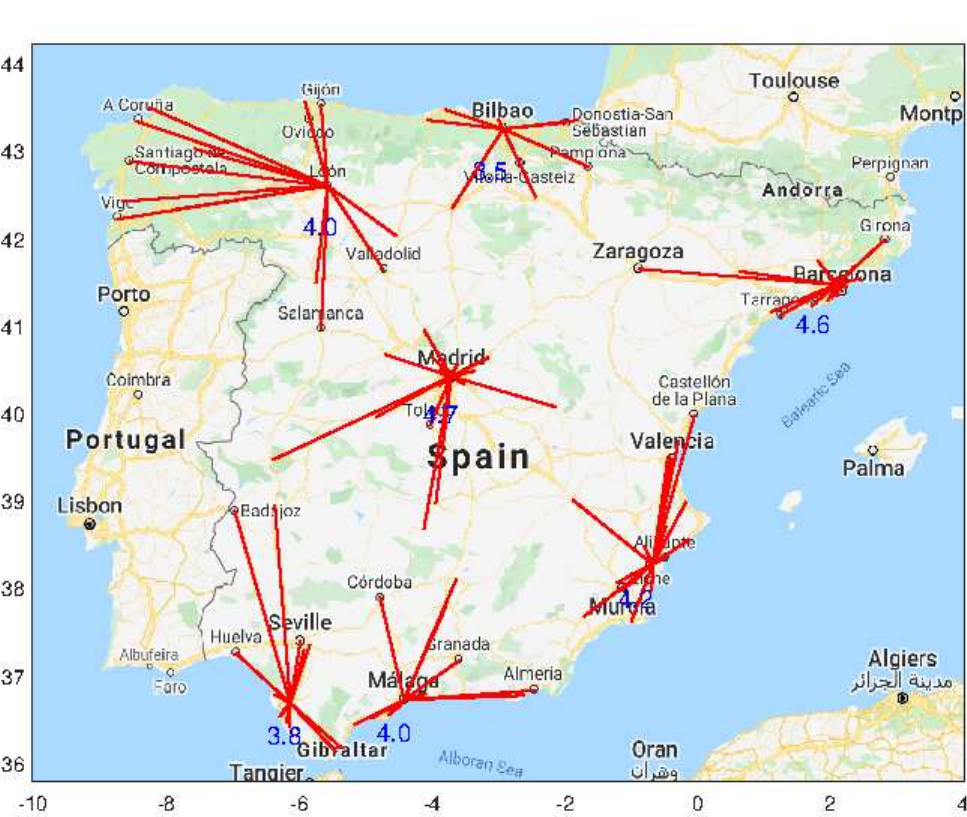}\\  (c)  Reduction to 7 facilities& (d) Reduction to 7 facilities (reallocation) \\
		\end{tabular}
		\caption{Results for 9 initial facilities without distribution costs. }\label{fig:Mapas}
	\end{center}
\end{figure}

Figure \ref{fig:Mapas} depicts the case without distributions costs. Panel (a) shows the original network and the nodes allocated to each facility. Panel (b) shows the results of the \textit{rebalancing} strategy; this strategy attains a 5.5\% improvement in profit and a reduction of 3.5\% in complexity costs. Two facilities were recentred, allowing the system to rebalance the demand allocation and, consequently, to reduce complexity in some facilities (Madrid and Barcelona specifically). The network resulting from the \textit{reduction} strategy is presented in panel (c). Two facilities were eliminated together with approximately 10\% of the total demand, bringing a 7.2\% reduction in complexity costs and a 4.5\% increase in profit. Finally, the network resulting after all demand uncovered nodes were reallocated to an open facility is shown in panel (d); this reallocation of nodes causes an increase in the cost of complexity of 5.1\%. 

\begin{figure}
	\centering
	\begin{tabular}{cc}
			\includegraphics[width=6cm]{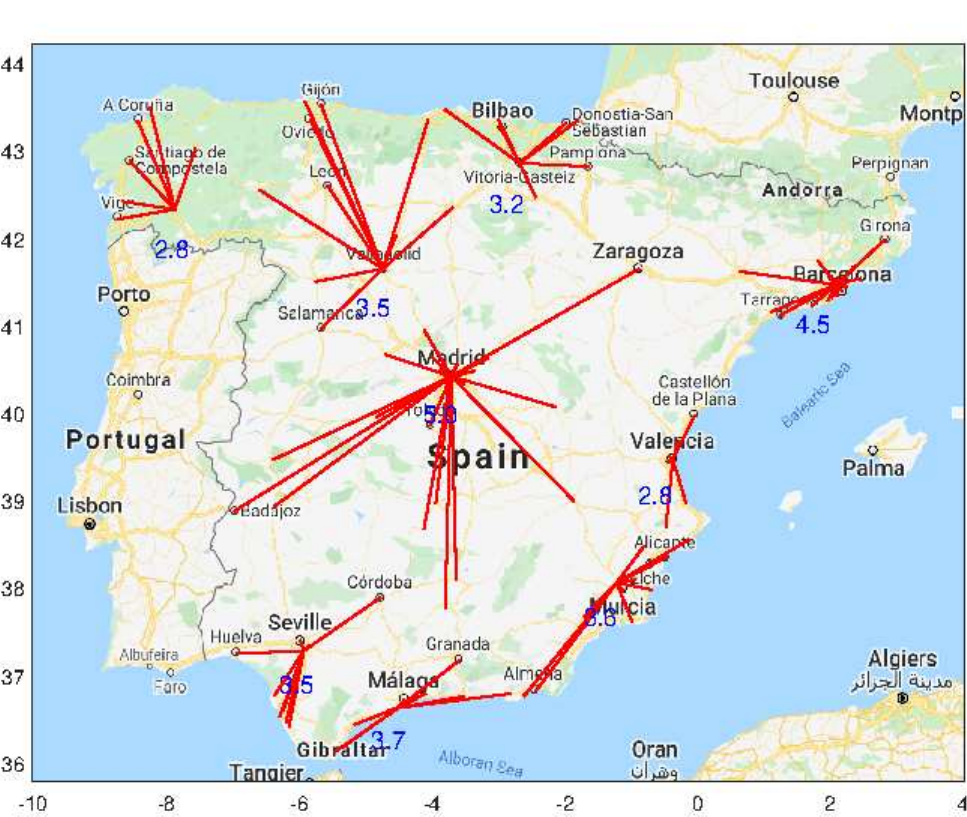} &
			\includegraphics[width=6cm]{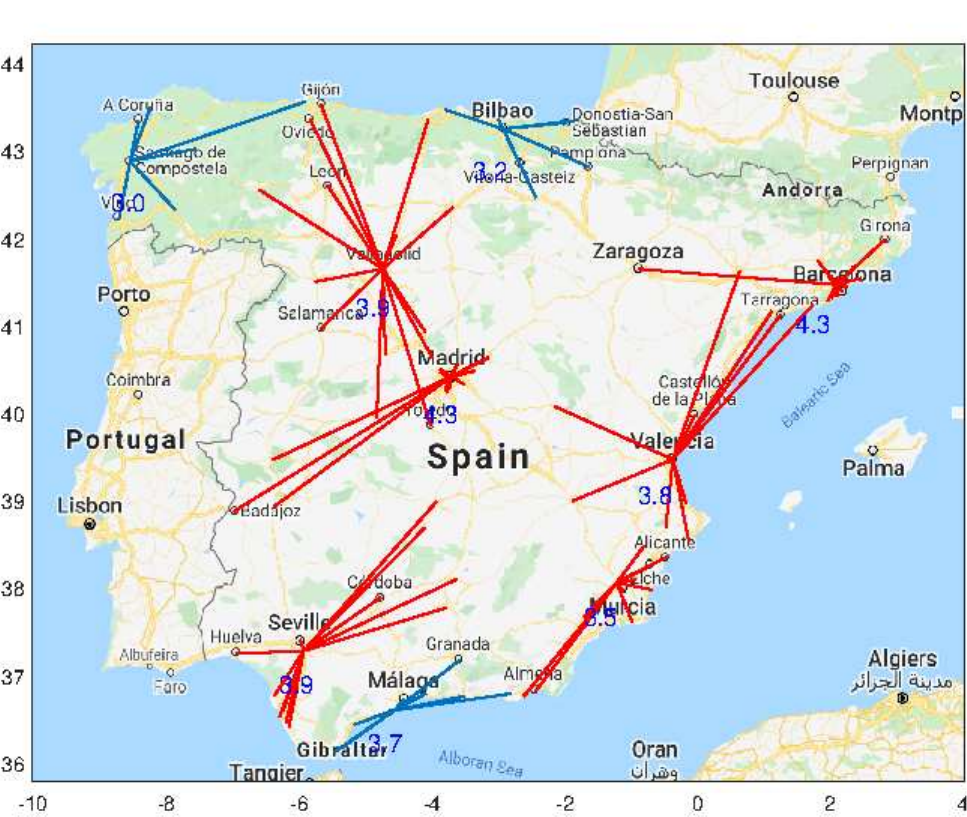} \\ (a) $Z^K$/$Z^K_{Plex}$ & (b) Rebalancing\\
			\includegraphics[width=6cm]{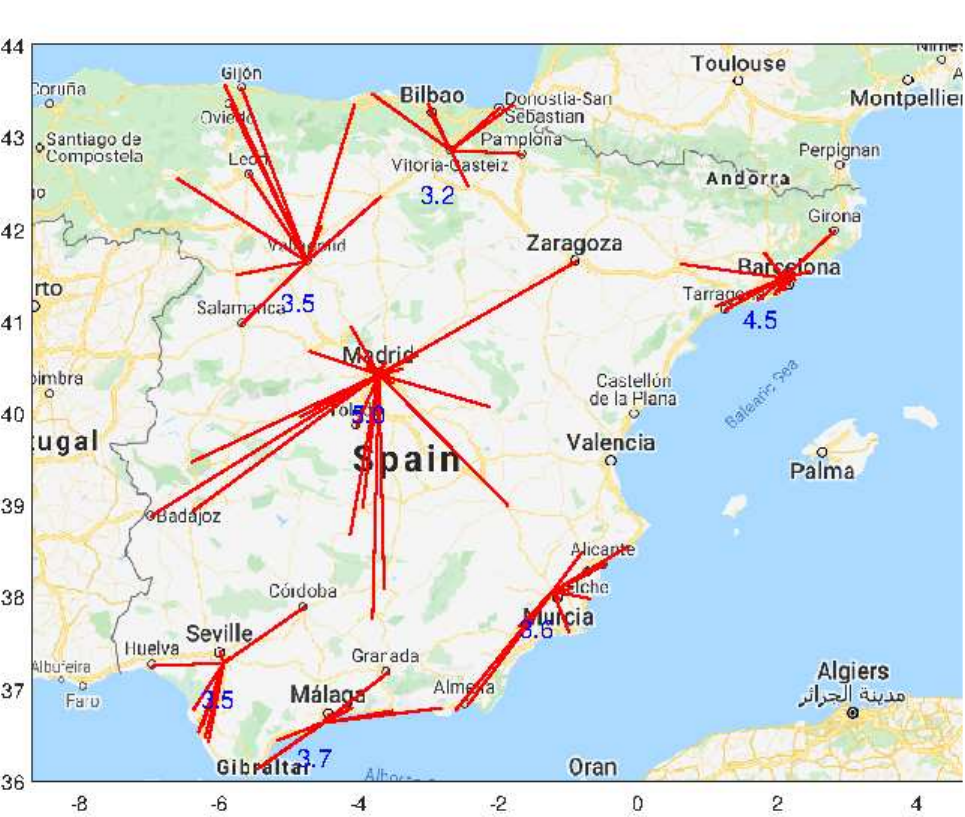}&
			\includegraphics[width=6cm]{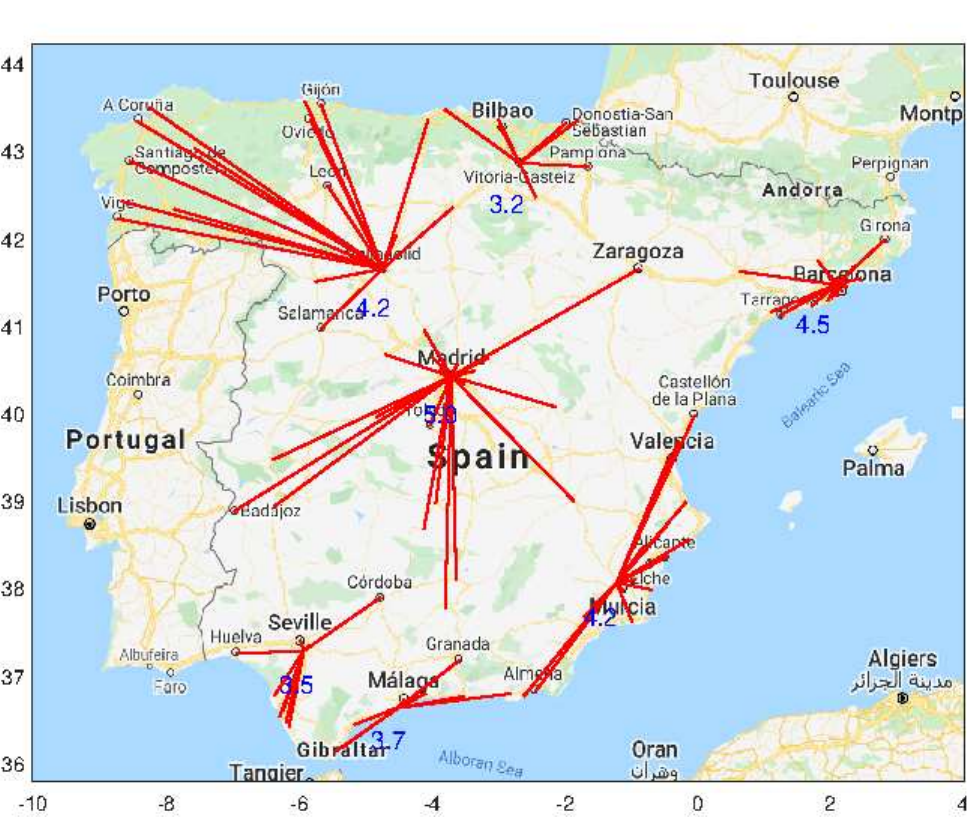}\\ (c)  Reduction to 7 facilities& (d) Reduction to 7 facilities (reallocation) \\
		\end{tabular}
		\caption{Results for 9 initial facilities with distribution costs. $\alpha = 0.125$, $\gamma= 200\ cts/km\times ton$, and $\varrho=66.6 \ cts/km\times ton$. }\label{fig:MapasDist}
\end{figure}

The results for the case with distribution costs are presented in Figure \ref{fig:MapasDist}. As before, panel (a) shows the original network. The map in panel (b) illustrated the resulting network after applying the \textit{rebalancing} strategy; in this case, three facilities were recentred, resulting in a reduction of 5.1\% in complexity costs. Panel (c) shows the network after applying the \textit{reduction} strategy; the same two facilities were eliminated, attaining a 7.7\% reduction in complexity costs and an 11\% increase in profit. Finally, panel (d) shows the network after uncovered demand nodes have been reallocated. The cost of complexity in this network is 5.3\% higher than in the initial network. 

\section{Conclusions}\label{sec:conc}

In this manuscript, we introduce the concept of \textit{locational complexity}, understood as the effect of an increasing number of facilities and their catchment areas on the company's performance. The main objective to explore the effectiveness of actions aimed at reducing location-related complexity in a distribution chain.

Using a measure for supply chain structural complexity, referred to as \textit{pars-complexity}, we develop a mathematical model that incorporates locational complexity and its costs in the decision process.

The mathematical properties of the \textit{pars}-complexity measure, allow us to establish that only by reducing both the number of facilities and the total demand covered it is possible to effectively reduce the network's locational complexity. Consequently, the impact of any network restructuring strategy which aims at reducing complexity-related costs without a significant reduction in market share, is limited to small improvements associated to a better balance of complexity-related costs among facilities.

The proposed mathematical model, is used to find an explanation to the problem of addiction to growth, and to explore the reasons why network restructuring strategies may result ineffective at reducing complexity-related costs. To accomplish this, we propose three alternative network restructuring strategies. In the first one, \textit{network rebalancing}, the company changes the allocation of demand across its facilities and, if feasible, relocates them seeking a more balanced distribution of distribution costs. In the second case, \textit{network rationalisation}, the company abandons a limited number of unprofitable demand nodes with the aim of reducing complexity associated costs, but the number of facilities remains the same. Finally, \textit{network reduction} strategy consists of eliminating a number of unprofitable facilities together with their associated demand. For completeness, we also analyse the case where demand from uncovered nodes is reallocated to facilities that remain open. Numerical experiments conducted on the three strategies, confirm our claim: network reduction strategies that are limited to eliminating facilities fail at reducing complexity-related costs.

Clearly, for many companies, losing market share is not always a viable or appealing option. Our experience suggests that a mechanism to reduce the total cost of complexity, without having to abandon entire markets, could be to divide the company into several, smaller, autonomous, and independent business units. This would guarantee that the (smaller) costs of complexity are absorbed by those smaller units. Modelling such a situation requires a deeper understanding about how location-related costs are generated in smaller units, and the development of optimisation and decision support models for this new problem is an avenue yet to be explored.

Finally, our findings highlight the importance of taking into account the network's structural complexity, and its costs, when developing network design models. Unfortunately, these models are non-linear and highly combinatorial, requiring the development of well-grounded heuristic approaches for their solution. This remains as an open line of research.

\newpage

\appendix

\section{Algorithms}

\subsection{The Network Rebalancing Algorithm}\label{alg:ZPlex}

\fbox{\parbox{12cm}{
	\begin{description}
	\item[Input:] $\mc S^K,\mathcal{N}^{\mc S^K}$
	\item[Do:] 
				\item  $\widetilde{Z}_{Plex}:=Z_{Plex}^\circ\left(\mathcal{N}^{\mc S^K}, \mc S^K\right)$
				\item[For all $\ell \in \mc S^K$] 
				\item $\mathcal{S}_\ell := \left\{ s_i \in \mathcal{N}_\ell: \tau\left(s_1\right) \geq \tau\left(s_2\right) \geq \ldots \geq \tau \left(s_{|\mathcal{N}_{\ell}|}\right)\right\}$
				\item[end for] 
				\item[For all] $\ell \in \mc S^K$ \\
				$\mathcal{M}:=\mathcal{N}^{\mc S^K}$; $i:=1$
				\begin{description}
					\item[While] $\nu=1$
					\begin{description}
						\item $\mathcal{N}'_\ell:=\mathcal{N}_\ell-\cur{s_i}$; \hspace{3pt} $\mathcal{N}'_{\varrho(s_i)}:=\mathcal{N}_{\varrho(s_i)} \cup s_i$
						\item $\overline{\mathcal{M}}:=\cur{\mathcal{M}-\left\{ \mathcal{N}_{\ell},\mathcal{N}_{\varrho(s_i)}\right\}} \cup \left\{ \mathcal{N}'_{\ell},\mathcal{N}'_{\varrho(s_i)}\right\}$ 
						\item $Z^{Imp}:=Z_{Plex}^\circ(\,\overline{\mc M}, \mc S^K)$
						\item[If] $Z^{Imp}\gg \widetilde{Z}_{Plex}$ 
						\begin{description}
							\item $\mathcal{M}:=\overline{\mathcal{M}}$ 
							\item $\mathcal{N}_{\ell}:=\mathcal{N}'_{\ell}$
							\item $\mathcal{N}_{\varrho(s_i)}:=\mathcal{N}'_{\varrho(s_i)}$
							\item $\widetilde{Z}_{Plex}:=Z^{Imp}$
							\item $i:=i+1$
						\end{description}
						\item[else]  $\nu:=0$ 
						\item[end if]
					\end{description}
					\item[end while]
				\end{description}
				\item[end for]
				\item[For all] $\ell \in \mc S^K$ \\
				$\pi\prn{\ell}=\underset{h \in \mathcal{N}_\ell }{\text{argmin}}\cur{\sum_{i \in \mathcal{N}_\ell} \gamma d(i,h)W_i}$
				\item[end for]
				\item $\overline{\mc S}:=\cur{\pi\prn{\ell}: \ell \in \mc S^K}$
				\item $\mathcal{N}^{\overline{\mc S}}:=\cur{\mathcal{N}_\ell: \ell \in \mc S^K}$
				\item $\widetilde{Z}_{Plex}:=Z_{Plex}^\circ\prn{\mathcal{N}^{\overline{\mc S}},\overline{\mc S}}$
				\item[Return:] $\widetilde{Z}_{Plex},\mathcal{N}^{\overline{\mc S}},\overline{\mc S}$
			\end{description}
}}

\vspace{6pt}

Some additional notation was required for this algorithm. For a given demand node $i \in \mathcal{N}_\ell$,
\[ \tau(i)=\prn{\min_{\ell' \in \mc S : \ell' \neq \ell}{\left\{d(i,\ell') \right\}}} W_i \]
provides the product of the distance from the node to its second closest facility times the node's demand; and
\[ \varrho(i)=\underset{\ell' \in \mc S : \ell' \neq \ell}{\text{argmin}}{\left\{d(i,\ell')  \right\}} \]
represents the second closest facility to node $i$.

Finally, for a given facility $\ell \in \mc S^K$ and its associated demand nodes $\mathcal{N}_\ell$, $\pi(\ell)$ represents the solution to the $1$-Median problem defined by $\ell$ and $\mathcal{N}_\ell$, i.e.
\[ \pi\prn{\ell}=\underset{h \in \mathcal{N}_\ell }{\text{argmin}}\cur{\sum_{i \in \mathcal{N}_\ell} \gamma d(i,h)W_i} \,. \]

\vspace{12pt}

\subsection{Network Rationalitation Algorithm}
\label{alg:DemRed}	

Let $\theta_\ell$ represet the ordered set of elements of $\mc N_\ell$, i.e. $\theta_\ell =\cur{\theta_{\ell_1},\theta_{\ell_2},\ldots,\theta_{\ell_{\abs{\mc N_\ell}}}}$, where $\theta_{\ell_k}=i\in \mc N_\ell : d(i,\ell)\geq d(\theta_{\ell_{k'}},\ell)$ for $k'=k+1,\ldots,\abs{\mc N_\ell}$. We can now define
$\Theta_\ell=\cur{\theta_{\ell_k} ,k=1,\ldots,n}$ as the set of the $n$ most distant nodes from facility $\ell$.

\vspace{6pt}

\fbox{\parbox{12cm}{

			\begin{description}
				\item[Input:] $\mc S^K,\mathcal{N}^{\mc S^K}$
				\item[For all] $\ell \in \mc S^K$
				\begin{description}
					\item Compute $\Theta_\ell$
					\item $\underline{\lambda}^0:=\left\{\lambda^0(i):i \in \Theta_\ell \right\}$
					\item $\Lambda^0_\ell := \text{argmax}\left\{ \lambda^0(i) : i \in \Theta_\ell \right\} $
					\item $\overline{\rule{0pt}{6pt}\mathcal{N}}:=\mathcal{N}^{\mc S^K}$
					\item $\overline{ \Theta}_\ell=\Theta_\ell$
					\item $\nu:=0$
					\item[While] $\Lambda^\nu_\ell >0$ $\vee$ $\overline{\Theta}_\ell \neq \emptyset$ 
					\begin{description}
						\item $\overline{\rule{0pt}{6pt}\mathcal{N}}:=\overline{\rule{0pt}{6pt}\mathcal{N}}-\left\{\Lambda^\nu_\ell \right\}$
						\item $\widetilde{Z}_{Plex}:= Z_{Plex}^\circ\left(\overline{\rule{0pt}{6pt}\mathcal{N}},\mc S^K \right) $
						\item $\overline{\Theta}_\ell:=\overline{\Theta}_\ell-{\Lambda^{\nu}_\ell}$
						\item $\nu:=\nu+1$
						\item $\underline{\lambda}^\nu:=\left\{\lambda^\nu(i):i \in \overline{\Theta}_\ell \right\}$
						\item $\Lambda^\nu_\ell := \text{argmax}\left\{ \lambda^\nu(i) : i \in \overline{\Theta}_\ell \right\}$
					\end{description}
					\item[end while]
				\end{description}
				\item[end for]
				\item[Return:] $\overline{\rule{0pt}{6pt}\mathcal{N}},\widetilde{Z}$
			\end{description}
    }}

\subsection{The Network Reduction Algorithm}\label{alg:NetRed}
\fbox{\parbox{12cm}{
	\begin{description}
				\item[Input:] $\mc S^K,\mathcal{N}^{\mc S^K}$
				\item[Do:] 
				\item  $\widetilde{Z}_{Plex}:=Z_{Plex}^\circ\left(\mathcal{N}^{\mc S^K}, \mc S^K\right)$
				\item $\mathcal{S}:=\mc S^K$
				\item $\mathcal{M}:=\mathcal{N}^{\mc S^K}$
				\item[While] $\nu = 1$
				\begin{description}
				\item[For all] $\ell \in \mathcal{S}$ \\
				$\mathcal{S}':=\mathcal{S}-\ell$ \\
				$\mathcal{M}':=\mathcal{M}-\mathcal{N}_\ell$ \\
				$Z_\ell=Z_{Plex}^\circ\left(\mathcal{M}', \mathcal{S}'\right)$ 
				\item[end for]
				\item $\ell^\circ  = \text{argmax}_{\ell}\cur{Z_\ell}$ 
				\item $Z^\circ = Z_{\ell^\circ}$ 
				\item[If] $Z^\circ > \widetilde Z_{Plex}$ \\
				$\widetilde Z_{Plex} := Z^\circ$ \\
				$\mathcal{S}:=\mathcal{S}-\ell^\circ$ \\
				$\mathcal{M}:= ,\mathcal{M}-\mathcal{N}_{\ell^\circ}$ 
				\item[else]
				$\nu:=0$
				\item[end if] 
                \end{description}
			\item[end while]
			\item[Return:] $\mathcal{S},\mathcal{M},\widetilde{Z}_{Plex}$
			\end{description}
}}

\bibliographystyle{unsrtnat}
\bibliography{ComplexityLoc}

\end{document}